\documentclass[11pt]{amsart}
\usepackage[top=1.3in, left=1.3in, right=1.3in, bottom=1.3in]{geometry}

\usepackage{amsmath,amsfonts,amsthm,graphicx}
\usepackage{amsmath}
\usepackage{hyperref}
\usepackage{graphicx,psfrag}  
\usepackage[psamsfonts]{amssymb}
\usepackage{amscd}
\usepackage[all,cmtip]{xy}
\title{What are Lyapunov exponents, and why are they interesting?}
\author{Amie Wilkinson}

\date{\today}
\theoremstyle{plain}

\newtheorem{theorem}{Theorem}[section]

\newtheorem{corollary}[theorem]{Corollary}

\newtheorem{dfinition}[theorem]{Definition}

\def\Diff{\operatorname{Diff} }

\def\mod{\hbox{mod} }

\def\title{\em}

\def\bar{\overline}

\def\cW{\mathcal{W}}

\def\cB{\mathcal{B}}

\def\cF{\mathcal{F}}

\def\cM{\mathcal{M}}

\def\C{\mathcal{C}}
\def\cC{\mathcal{C}}

\def\cZ{\mathcal{Z}}

\def\dist{\operatorname{dist}}
\def\transverse{\,\raise2pt\hbox to1em{\hfil$\top$\hfil}\hskip -1em \hbox
to1em{\hfil$\cap$\hfil}\,} 

\newcommand\RR{{\mathbb R}}
\newcommand\CC{{\mathbb C}}

\newcommand\HH{{\mathbb H}}

\newcommand\TT{{\mathbb T}}
\newcommand\NN{{\mathbb N}}
\newcommand\ZZ{{\mathbb Z}}

\newlength{\figboxwidth} 

\begin{document}
\maketitle

\section*{Introduction}

At the 2014 International Congress of Mathematicians in Seoul, South Korea, Franco-Brazilian mathematician Artur Avila was awarded the Fields Medal for ``his profound contributions to dynamical systems theory, which have changed the face of the field, using the powerful idea of renormalization as a unifying principle."\footnote{\url{http://www.mathunion.org/general/prizes/2014/prize-citations/}}    Although it is not explicitly mentioned in this citation, there is a second unifying concept in Avila's work that is closely tied with renormalization: Lyapunov (or characteristic) exponents.  Lyapunov exponents play a key role in three areas of Avila's research: smooth ergodic theory,
billiards and translation surfaces, and the spectral theory of 1-dimensional Schr\"odinger operators.  Here we take the opportunity to explore these areas and reveal some underlying themes connecting exponents, chaotic dynamics and renormalization.

But first, what are Lyapunov exponents?  Let's begin by viewing them in one of their natural habitats: the iterated barycentric subdivision of a triangle.   

When the midpoint of each side of a triangle is connected to its opposite vertex by a line segment, the three resulting segments meet in a point in the interior of the triangle.  The barycentric subdivision of a triangle is the collection of 6 smaller triangles determined by these segments and the edges of the original triangle:
\begin{figure}[ht]
\begin{center}
\includegraphics[scale=.4 ]{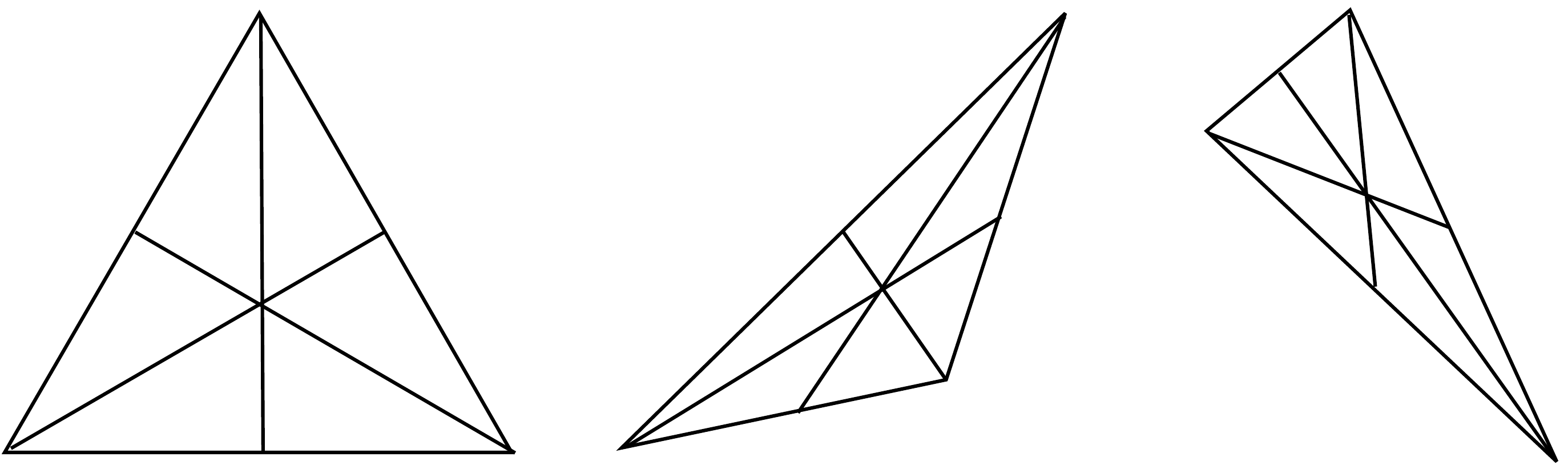}
\end{center}
\caption{Barycentric subdivision.}
\end{figure}

As the process of barycentric subdivision starts with a triangle and produces triangles, it's natural to iterate the process, barycentrically subdividing the six triangles obtained at the first step, obtaining 36 triangles, and so on, as in Figure 2.
\begin{figure}[ht]
\begin{center}
\includegraphics[scale=.4 ]{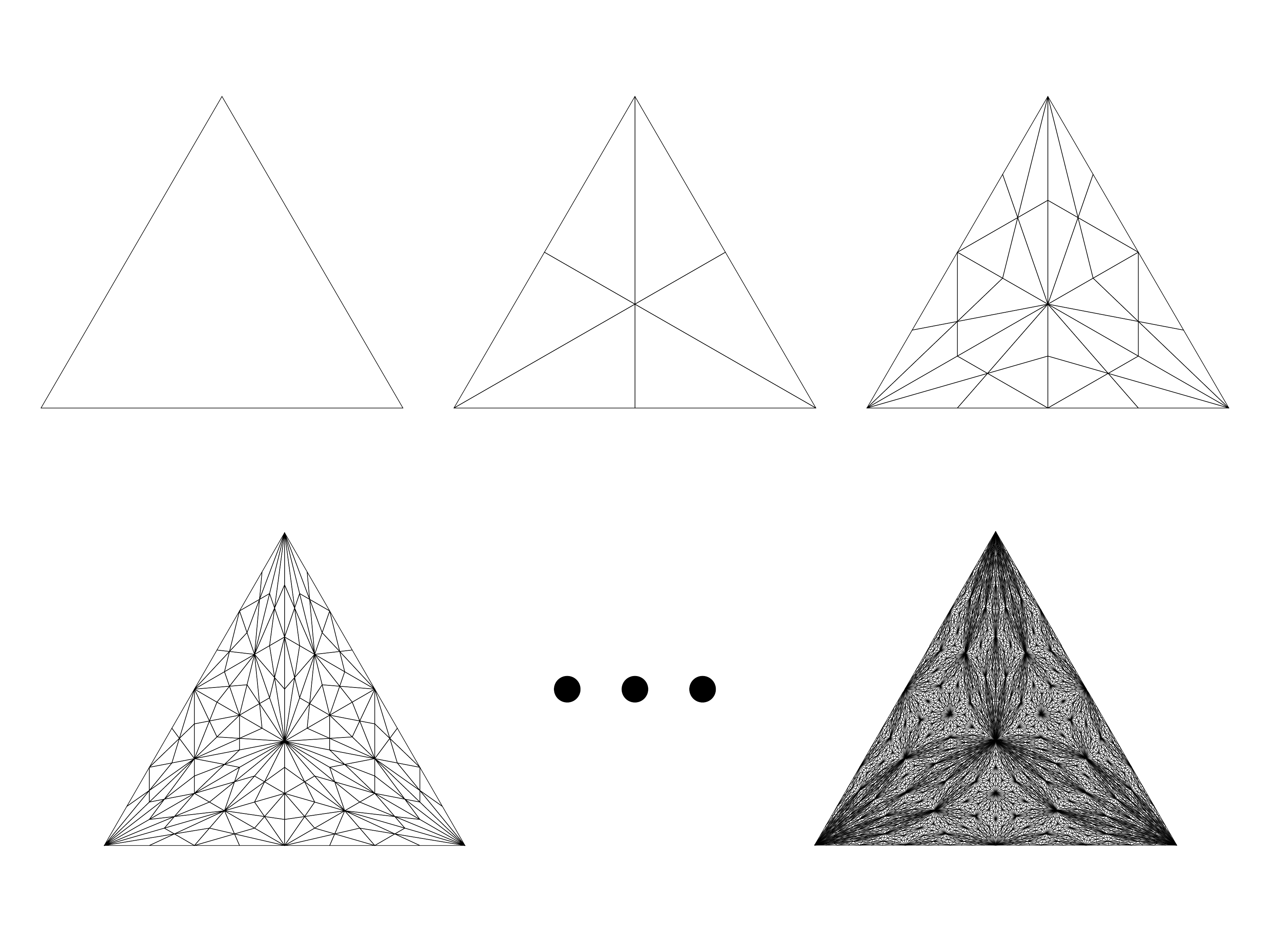}
\end{center}
\caption{Iterating barycentric subdivision, from \cite{McMo}.}
\end{figure}

Notice that as the subdivision gets successively finer,  many of the triangles produced by subdivision get increasingly eccentric and needle-like.  We can measure the skinniness of a triangle $T$ via the {\em aspect ratio} $\alpha(T) = \hbox{area}(T)/L(T)^2$, where $L(T)$ is the maximum of the side lengths; observe that similar triangles have the same aspect ratio. Suppose we  fix a rule for labeling the triangles in a possible subdivision $1$ through $6$, roll a six-sided fair die and at each stage choose a triangle to subdivide.  The sequence of triangles $T_1 \supset T_2 \supset  \ldots$ obtained have aspect ratios $\alpha_1, \alpha_2, \ldots$, where $\alpha_n=\alpha(T_n)$.
\begin{theorem}[\cite{McMo}, see also \cite{Barany:Beardon:Carne}] There exists a real number $\chi \approx  0.0446945 >0$ such that almost surely: \[\lim_{n\to\infty} \frac{1}{n}\log\alpha_n =  - 2 \chi.\]
\end{theorem}
In other words,  if triangles are chosen successively by a random  toss of the die, then with probability 1, their aspect ratios will tend to $0$ in the $n$th toss at 
an exponential rate governed by $\exp(-2 n\chi)$.  The same conclusion holds with the same value of $\chi$ if the initial equilateral triangle is replaced by any marked triangle.  This magical number $\chi$ is a Lyapunov exponent.  We return to this example at the end of the next section. 

Lyapunov exponents make multiple appearances in the analysis of dynamical systems.   After defining basic concepts and explaining examples in Section 1, we describe in Sections 2--4 a sampling of
Avila's results in smooth ergodic theory, Teichm\"uller theory and spectral theory, all of them tied to Lyapunov exponents in a fundamental way.  We explore some commonalities of these results in Section  5.  Section 6 is  devoted to a discussion of some themes that arise in connection with Lyapunov exponents.

\section{Cocycles, exponents and hyperbolicity}

Formally, Lyapunov exponents are quantities associated to a cocycle over a measure-preserving dynamical system.  A measure-preserving dynamical system is a triple $(\Omega, \mu, f)$, where  $(\Omega,\mu)$  is a probability space, and  $f\colon \Omega\to \Omega$ is a measurable map that preserves the measure $\mu$, meaning that $\mu(f^{-1}(X))=\mu(X)$ for every measurable set $X\subset \Omega$.  We say that $(\Omega,\mu,f)$ is {\em ergodic} if the only $f$-invariant measurable sets have  $\mu$-measure $0$ or $1$.  Equivalently,
 $(\Omega,\mu,f)$ is ergodic if  the only functions $\phi\in L^2(\Omega,\mu)$ satisfying  $\phi\circ  f = \phi$ are  the constant functions.
 Any $f$-invariant measure $\mu$ can be canonically decomposed into ergodic invariant measures, a fact that allows us to restrict attention to ergodic measures in some contexts, simplifying statements.  The measures in such a decomposition are called {\em ergodic components,} and there can be uncountably many of them.  The process of ergodic decomposition is a bit technical to describe; we refer the reader to \cite{Manebook} for more details.

\subsection{Examples of measure-preserving systems}

Here is a list of  three types of measure-preserving systems that we will refer to again in the sections that follow.

\medskip

\noindent  {\bf Rotations on the circle.}  On the circle $\Omega = \RR/\ZZ$,
let $f_\alpha(x)= x+\alpha \, (\mod\, 1)$, where $\alpha\in \RR$ is fixed.  The map $f_\alpha$ preserves the Lebesgue-Haar measure  $\mu$ (that assigns to an interval $I\in \RR/\ZZ$ its length $|I|$). The map $f_\alpha$ is ergodic with respect to $\mu$ if and only if $\alpha$ is irrational. This has a straightforward proof: consider the equation
$\phi\circ f_\alpha = \phi$, for some $\phi\in L^2(\RR/\ZZ,\mu)$,  and solve for the Fourier coefficients of $\phi$.  

When $\alpha = p/q$ is rational, every point $\omega\in \Omega$ is periodic,  satisfying $f^q(\omega) = \omega$.  Each  $\omega$ then determines an ergodic $f_\alpha$-invariant probability measure $\nu_\omega$ obtained by averaging the Dirac masses along the orbit of $\omega$:
\[
\nu_\omega := \frac{1}{q}\left(\delta_\omega + \delta_{f_\alpha(\omega)} + \cdots + \delta_{f_\alpha^{q-1}(\omega)} \right).
\]
Each $\nu_\omega$ is an ergodic component of the measure $\mu$, and hence there are uncountably many such components.

When $\alpha$ is irrational, $\mu$ is the {\em unique}  $f_\alpha$-invariant Borel probability measure. A homeomorphism of a compact metric space that has a unique invariant measure is called {\em uniquely ergodic} --- a property that implies ergodicity and more.  Unique ergodicity is mentioned again in Section~\ref{s=lyaperg} and is especially relevant to the discussion of Schr\"odinger operators with quasiperiodic potentials in Section~\ref{s=Schrodinger}.

There is nothing particularly special about the circle, and the properties of circle rotations listed here generalize easily to rotations on compact abelian groups.

\medskip

\noindent {\bf Toral automorphisms.}  Let $\Omega = \TT^2 := \RR^2/\ZZ^2$, the $2$-torus. Fix a matrix $A\in SL(2,\ZZ)$. Then $A$ acts linearly on the plane by multiplication and preserves the lattice $\ZZ^2$, by virtue of having integer entries and determinant $1$.  It therefore induces 
a map  $f_A\colon \TT^2\to  \TT^2$ of the 2-torus, a group automorphism.
The area $\mu$ is preserved because $\det(A)=1$.  Such an automorphism is ergodic with respect to $\mu$ if and only if the eigenvalues of $A$ are $\lambda$ and $\lambda^{-1}$, with
$|\lambda| > 1$.
This can be proved by examining the Fourier coefficients of $\phi\in L^2(\TT^2, \mu)$ satisfying $\phi\circ f_A = \phi$: composing with $f_A$ permutes the Fourier coefficients of $\phi$ by the adjoint action of $A$, and the assumption on $A$ implies that if $\phi$ is not constant there must be infinitely many identical coefficients, which violates the assumption that $\phi\in L^2(\TT^2, \mu)$.\footnote{In higher dimensions, a matrix $A\in SL(d,\ZZ)$ similarly induces an automorphism $f_A$ of $\TT^n := \RR^n/\ZZ^n$.  The same argument using Fourier series shows that $f_A$ is ergodic if and only if $A$ does not have a root of $1$ as an eigenvalue.}  

 In contrast to the irrational rotation $f_\alpha$, the map $f_A$ has many invariant Borel probability measures, even when $f_A$ is ergodic with respect to the area $\mu$.  For example, as we have just seen, every periodic point of $f_A$ determines an ergodic invariant measure, and $f_A$ has infinitely many periodic points.
This is a simple consequence of the Pigeonhole Principle, using the fact that $A\in SL(2,\ZZ)$:
for every natural number $q$, the finite collection of points $\{(\frac{p_1}{q}, \frac{p_2}{q}): p_1, p_2\in \{0,\ldots, q-1\}\}\subset \TT^2$ is permuted by $f_A$, and so each element of this set is fixed by some power of $f_A$.

\medskip

\noindent{\bf Bernoulli shifts.}  Let $\Omega= \{1,\ldots, k\}^\NN$ be the set of all infinite, one sided strings $\omega = (\omega_1,\omega_2,\cdots )$  on the alphabet  $\{1,\ldots, k\}$. Endowed with the product topology, the space $\Omega$ is compact, metrizable, and homeomorphic to a Cantor set.  The shift map $\sigma\colon \Omega\to\Omega$
is defined by $\sigma(\omega)_k = \omega_{k+1}$.  In other words, the image of the sequence
$\omega = (\omega_1,\omega_2,\cdots )$   is the shifted sequence  $\sigma(\omega)= (\omega_2,\omega_3,\cdots )$.  Any nontrivial probability vector $p=(p_1,\ldots, p_k)$  (i.e. with $p_i\in (0,1)$, and $\sum_i p_i = 1$) defines a product measure
$\mu = p^\NN$  supported on $\Omega$.\footnote{The product measure has a simple description in this context: setting $C_i(j):=  \{\omega\in \Omega: \omega_i = j\}$,  the measure $\mu = (p_1,\ldots, p_k)^\NN$ is defined by the properties:  $\mu(C_i(j)) = p_j$, and
 $\mu(C_i(j) \cap C_{i'}(j')) = p_j p_{j'}$, for any   $i\neq i' \in \NN$ and $j, j' \in\{1,\ldots,k\}$. }  The triple $\left(\Omega, \sigma, \mu = (p_1,\ldots, p_k)^\NN\right)$ is called a {\em Bernoulli shift},
and $\mu$ is called a {\em Bernoulli measure}.  It is not hard to see that the shift $\sigma$ preserves $\mu$ and is ergodic with respect to $\mu$.

The shift map $\sigma\colon \Omega\to \Omega$  manifestly has uncountably many 
invariant Borel probability measures, in particular the Bernoulli measures $(p_1,\ldots, p_k)^\NN$, but the list does not end there.  In addition to periodic measures (supported on orbits of periodic strings $(\omega_1, \ldots, \omega_q, \omega_1, \ldots, \omega_q,  \ldots)$), there are $\sigma$-invariant  probability measures on $\Omega$ encoding every measure preserving continuous dynamical system\footnote{Subject to some constraints involving invertibility and entropy.} --- in this sense the shift  is a type of universal dynamical system.

\subsection{Cocycles}\label{ss=cocycles}

Let $M_{d\times d}$ be the $d^2$-dimensional vector space of $d\times d$ matrices (real or complex).  A {\em cocycle} is a pair $(f,A)$, where $f\colon \Omega\to \Omega$ and  $A\colon \Omega\to M_{d\times d}$ are measurable maps. We also say that $A$  is a {\em cocycle over $f$}.  For each $n>0$, and $\omega\in\Omega$, we write
$$A^{(n)}(\omega) = A(f^{n-1}(\omega))A(f^{n-2}(\omega))\cdots A(f(\omega)) A(\omega),$$
where $f^n$ denotes the $n$-fold composition of $f$ with itself.
For $n=0$, we set $A^{(n)}(\omega) = I$, and 
if both $f$ and the values of the cocycle $A$ are invertible, we also define, for $n\geq 1$:
$$
A^{(-n)}(\omega) =  (A^{(n)}(f^{-n+1}(\omega)))^{-1} = A^{-1}(f^{-n+1}(\omega))\cdots A^{-1}(\omega).
$$

One comment about the terminology ``cocycle:"  while $A$ is colloquially referred to as a cocycle over $f$, to
fit this definition into a proper cohomology theory, 
one should reserve the term cocycle for the function
$(n, \omega)\to A^{(n)}(\omega)$  and call $A$ the {\em generator} of this (1-)cocycle.  See \cite{ASVW} for a more thorough discussion of  this point.

A fruitful way of viewing a cocycle $A$ over $f$ is as a hybrid dynamical system  $(f,A) \colon \Omega\times M_{d\times d} \to  \Omega\times M_{d\times d}$ defined by
\[ 
(f,A) (\omega, B) = (f(\omega), A(\omega) B).
\]
Note that the $n$th iterate $(f,A)^n$ of this hybrid map is the hybrid map $(f^n, A^{(n)})$.
The vector bundle $\Omega\times M_{d\times d}$ can be reduced in various ways to obtain associated hybrid systems, for example, the map $(f,A)\colon \Omega\times \RR^d \to \Omega\times \RR^d$ defined by 
$(f,A) (\omega, v) = (f(\omega), A(\omega) v)$.   Thus a natural generalization  of a cocycle over $f$ is 
a map $F\colon \cB\to \cB$, where $\pi\colon \cB\to \Omega$ is a vector bundle, and $F$ acts linearly on fibers, with $\pi\circ F = f\circ \pi$.  We will use this extended definition of cocycle
to define the  derivative cocycle in Subsection~\ref{ss=examples}.

\subsection{Lyapunov exponents}

Let $f\colon \Omega\to \Omega$ be a measurable map (not necessarily preserving a probability measure).
We say that a real number $\chi$ is a {\em  Lyapunov exponent for the cocycle $A\colon \Omega\to M_{d\times d}$ over $f$
at the point $\omega\in \Omega$}  if there exists a nonzero vector  $v\in\RR^d$, such that
\begin{equation}\label{e=lyapdef}
\lim_{n\to \infty} \frac{1}{n} \log\|A^{(n)}(\omega) v\| = \chi.
\end{equation}
Here $\|\cdot\|$ is a fixed norm on the vector space space $M_{d\times d}$.  The limit in (\ref{e=lyapdef}), when it exists,  does not depend on the choice of such a norm (exercise).

Oseledets proved in 1968 \cite{Os} that  if $(f,\Omega,\mu)$ is a measure-preserving system, and
 $A$ is a cocycle over $f$ satisfying  the integrability condition $\log\|A\|  \in L^1(\Omega,\mu)$, 
then for $\mu$-almost every $\omega\in \Omega$ and for every nonzero $v\in \RR^d$ the limit in (\ref{e=lyapdef}) exists.  This limit assumes at most $d$ distinct values $\chi_1(\omega) > \chi_2(\omega) > \cdots > \chi_{k(\omega)}(\omega)$.
Each exponent $\chi_i(\omega)$ is achieved with a multiplicity $d_i(\omega)$ equal to the dimension of the space of vectors $v$ satisying (\ref{e=lyapdef}) with $\chi = \chi_i$, 
and these multiplicities satisfy $\sum_{i=1}^{k(\omega)} d_i(\omega) = d$.

 If the cocycle $A$ takes values in $SL(d,\RR)$,  then, since $\log\det(A(\omega)) \equiv \log(1) = 0$,
we obtain that $\sum_{i=1}^{k(\omega)} d_i(\omega)\chi_i(\omega)= 0$. Thus if $A$ takes values in $SL(2,\RR)$, then the exponents are of the form $-\chi(\omega)\leq 0 \leq \chi(\omega)$. 

If $(\Omega,\mu,f)$ is ergodic, then the functions $k(\omega)$, $\chi_i(\omega)$ and $d_i(\omega)$,  are constant $\mu$-almost everywhere.  In this case, the essential values  $\chi_1,\dots, \chi_k\in \RR$ are called the {\em Lyapunov exponents of $A$ with respect to the ergodic measure $\mu$}.

\subsection{Two important classes of cocycles}\label{ss=examples}

{\em Random matrix cocycles}  encode the behavior of a random product of matrices.
Let $\{A_1, \ldots, A_k\} \subset  M_{d\times d}$ be a finite collection of matrices.
Suppose we take a $k$-sided die and roll it repeatedly.  If the
\begin{figure}[ht]
\begin{center}
\includegraphics[scale=0.23]{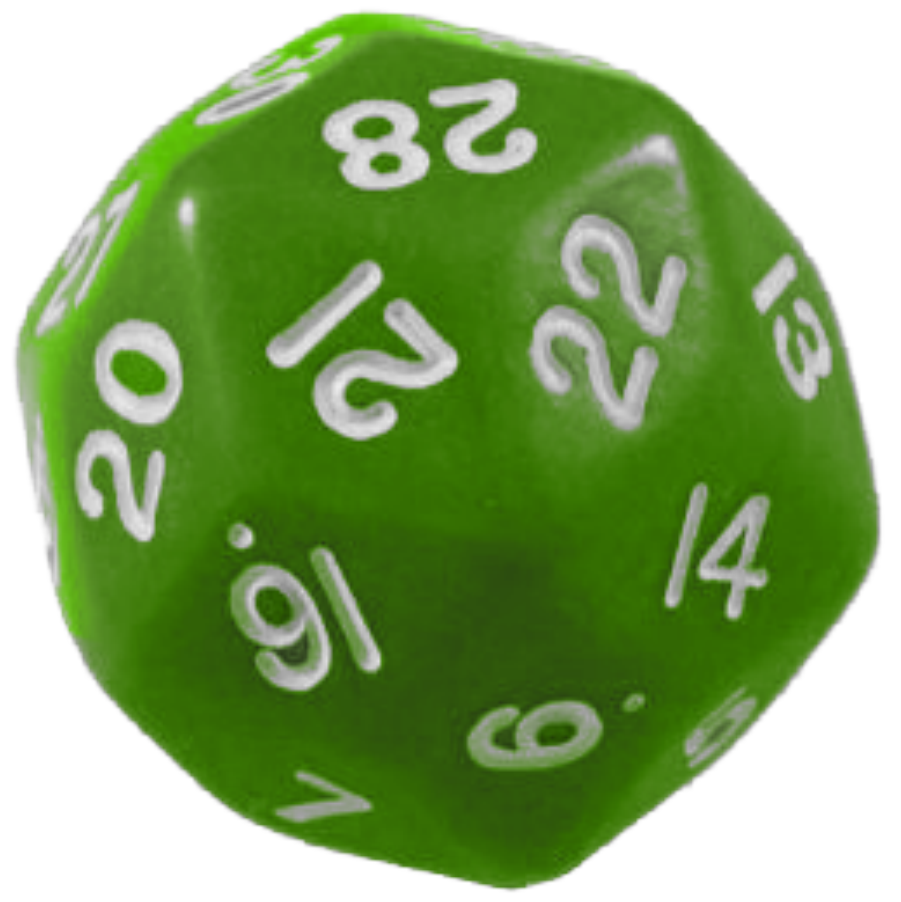}
\end{center}
\caption{A 30-sided Dungeons and Dragons die.}
\end{figure}
die comes up with the number $j$, we choose the matrix $A_j$, thus creating a sequence
$A_{\omega_1}, A_{\omega_2},\ldots$, where $\omega = (\omega_1,\omega_2,\ldots)\in \{1, \ldots,k\}^\NN$. This process can be packaged in a cocycle $A$ over a measure preserving system $(\Omega, \mu, \sigma)$ by setting $\Omega = \{1,\ldots, k\}^\NN$, $\mu = (p_1,\ldots, p_k)^\NN$,  where $p_j$ is the probability that the die shows
$j$ on a roll, and  setting $\sigma$ to be the shift map.
The cocycle is  defined by $A(\omega) = A_{\omega_1}$.
Then $A^{(n)}(\omega)$ is simply the product of the first $n$ matrices produced by this process.

More generally, suppose that
 $\eta$ is a probability measure on the set of matrices $M_{d\times d}$. The 
space $\Omega = M_{d\times d}^\NN$ of sequences $(M_1, M_2,  \dots)$ carries the product measure $\eta^\NN$, which is invariant under the shift map $\sigma$, where as above $\sigma(M_1, M_2,  \dots) = (M_2, M_3,  \dots)$.  There is a natural cocycle $A\colon \Omega\to  M_{d\times d}$ given by $A((M_1, M_2,  \dots)) = M_1$.  The matrices
$A^{(n)}(\omega)$, for $\omega\in \Omega$ are just $n$-fold random products of matrices chosen independently with respect to the measure $\eta$.

In the study of smooth dynamical systems the {\em derivative cocycle} is a central player. Let  $f\colon M\to M$ be a $C^1$ map on a compact $d$-manifold $M$.  Suppose for simplicity that the tangent bundle is trivial:
$TM = M\times \RR^d$.
Then for each $x\in M$, the derivative $D_xf\colon T_xM\to T_{fx} M$  can be written as a matrix $D_xf \in M_{d\times d}$. The map $x\mapsto D_xf$  is called the derivative cocycle. 
The Chain Rule implies that if $A=Df$ is 
a derivative cocycle, then
$D_xf^n = A^{(n)}(x)$.  

The case where $TM$ is not trivializable is easily treated: either one trivializes $TM$ over a suitable subset of $M$, or one expands the definition of cocycle as described  at the end of Subsection~\ref{ss=cocycles}: the map
$Df\colon TM\to TM$ is an automorphism of the vector bundle $TM$, covering the map $f$.  Lyapunov exponents for the derivative cocycle are defined analogously to (\ref{e=lyapdef}).  We fix  a continuous choice of norms
$\{\|\cdot\|_ x\colon T_xM \to \RR_{\geq 0}: x\in M\}$, for example the norms given by a Riemannian metric (more generally, such a family of norms is called a {\em Finsler}).  Then $\chi$ is a Lyapunov exponent for $Df$ at $x\in M$ if there exists $v\in T_xM$ such that
\begin{equation}\label{e=lyapdefderiv}
\lim_{n\to \infty} \frac{1}{n} \log\|D_xf^n  v\|_{f^n(x)} = \chi.
\end{equation}
Since $M$ is compact, the Lyapunov exponents of $Df$ do not depend on the choice of Finsler.  The conclusions of Oseledets's theorem hold analogously for derivative cocycles with respect to any $f$-invariant measure on $M$.

A simple example of a derivative cocycle is provided by the toral automorphism $f_A\colon \TT^2\to \TT^2$ described above.  Conveniently, the tangent bundle to $\TT^2$ is trivial, and the derivative cocycle is the constant cocycle $D_xf_A = A$.

\subsection{Uniformly hyperbolic cocycles}

A  special class of cocycles whose Lyapunov exponents are nonzero with respect to {\em any} invariant probability measure
are the uniformly hyperbolic cocycles.
 
\begin{dfinition}\label{d=unifhyp} A continuous cocycle $A\colon \Omega\to M_{d\times d}$ over a homeomorphism $f\colon \Omega\to\Omega$
of a compact metric space $\Omega$ is {\em uniformly hyperbolic} if there exists an integer $n\geq 1$, and for every $\omega\in \Omega$, there is a splitting
$\RR^d = E^u(\omega)\oplus E^s(\omega)$ into subspaces that depend continuously on $\omega$, such that for every $\omega \in \Omega$:
\begin{itemize}
\item[(i)]  $A(\omega)E^u(\omega) = E^u(f(\omega))$, and  $A(\omega)E^s(\omega) = E^s(f(\omega))$, 
\item[(ii)]  $v\in E^u(\omega)\implies \|A^{(n)}(\omega)v\| \geq 2\|v\|$, and
\item[(iii)]  $v\in E^s(\omega)\implies \|A^{(-n)}(\omega) v\| \geq 2\|v\|$.
\end{itemize}
\end{dfinition}
The definition is independent of choice of norm $\|\cdot\|$; changing norm on $M_{d\times d}$ simply changes the value of $n$.  The number $2$ in conditions (ii) and (iii) can be replaced by any fixed real number greater than $1$; again this only changes the value of $n$.
Notice that measure plays no role in the definition of uniform hyperbolicity.  It is a topological
property of the cocycle.  For short, we say that $(f,A)$ is uniformly hyperbolic.

A trivial example of a uniformly hyperbolic cocycle is the constant cocycle $A(\omega) \equiv A_0\in SL(2,\RR)$,
where $A_0$ is any matrix whose eigenvalues $\lambda > \lambda^{-1}$ satisfy $\lambda > 1$.  Here the splitting 
$\RR^2 = E^u(\omega)\oplus E^s(\omega)$ is the constant splitting into the sum of the $\lambda$ and $\lambda^{-1}$ eigenspaces of $A_0$, respectively.  For a constant cocycle, the Lyapunov exponents
are defined everywhere and are also constant; for this $SL(2,\RR)$ cocycle, the exponents are $\pm\log\lambda$.

A nontrivial example of a uniformly hyperbolic cocycle is any nonconstant, continuous $A\colon \Omega\to SL(2,\RR)$ with the property that the entries of $A(\omega)$ are all positive, for any $\omega\in \Omega$.  In this case the splitting is given by
\[
E^u(\omega) := \bigcap_{n\geq 0} A^{(n)}(f^{-n}(\omega))\left(\C_+ \right), \;\hbox{ and }\; E^u(\omega) := \bigcap_{n\geq 0} A^{(-n)}(f^{n}(\omega))\left(\C_- \right),
\]
where $\C_+$ denotes the set of $(x,y)\in \RR^2$ with $xy\geq 0$, and  $\C_-$  is the set of $(x,y)$ with $xy\leq 0$.  
For an example of this type, the Lyapunov exponents might not be everywhere defined, and their exact values with respect to a particular invariant measure are not easily determined, although they will always be nonzero where they exist (exercise).
In this example and the previous one, the base dynamics  $f\colon \Omega\to \Omega$ are irrelevant as far as uniform hyperbolicity of the cocycle is concerned.

Hyperbolicity is an {\em open} property of both the cocycle $A$ and the dynamics $f$:  if  $(f,A)$ is uniformly hyperbolic,
and $\hat f$ and $\hat A$ and both uniformly close (i.e. in the $C^0$ metric)  to  $f$ and $A$,  then  $(\hat f, \hat A)$ is uniformly hyperbolic.  The reason is that, as in the example just presented, uniform hyperbolicity is equivalent to the existence of continuously varying cone families
$\cC_+(\omega), \cC_-(\omega) \subset \RR^d$, jointly spanning $\RR^d$ for each $\omega\in \Omega$, and an integer $n\geq 1$ with the following properties:  $A^{(n)}(\omega)\left(\cC_+(\omega)\right) \subset \cC_+(f^n(\omega))$;  $A^{(-n)}(\omega)\left(\cC_-(\omega)\right) \subset \cC_-(f^{-n}(\omega))$;   vectors in $\cC_+(\omega)$ are  doubled in length by $A^{(n)}(\omega)$; and vectors in  $\cC_-(\omega)$ are  doubled in length by $A^{(-n)}(\omega)$.
 The existence of such cone families is preserved under small perturbation.

\subsection{Anosov diffeomorphisms}\label{ss=Anosov}

A diffeomorphism $f\colon M\to M$  whose derivative cocycle is uniformly hyperbolic is called {\em Anosov}. Again, one needs to modify this definition when the tangent bundle $TM$ is nontrivial: the splitting of $\RR^d$ in the definition is replaced by a splitting $TM = E^u\oplus E^s$ into subbundles --- that is,  a splitting $T_xM = E^u(x) \oplus E^s(x)$ into subspaces, for each $x\in M$, depending continuously on $x$. The norm $\|\cdot\|$ on the space $M_{d\times d}$  is replaced by a Finsler, as in the discussion at the end of Subsection~\ref{ss=examples}.   Since $M$ is assumed to be compact,  the Anosov property does not depend on the choice of Finsler.

Anosov diffeomorphisms remain Anosov after a  $C^1$-small perturbation, by the openness of uniform hyperbolicity  of cocycles.   More precisely, the $C^1$ distance $d_{C^1}(f,g)$ between two diffeomorphisms is the sum of the  $C^0$ distance between $f$ and $g$ and the $C^0$ distance between $Df$ and $D g$; thus if $f$ is Anosov
and $d_{C^1}(f,g)$ is sufficiently small, then $Dg$ is hyperbolic, and so $g$ is Anosov.  Such a $g$ is often called a {\em $C^1$ small perturbation of $f$.}

The toral automorphism $f_A\colon \TT^2\to\TT^2$,
with $A=\begin{pmatrix} 2&1\\1&1\end{pmatrix}$ is Anosov; since the derivative cocycle is constant, the splitting $\RR^2 = E^u(x) \oplus E^s(x)$, for $x\in \TT^2$ does not depend on $x$: as above, $E^u(x)$ is the expanding eigenspace for $A$ corresponding to the larger eigenvalue $\lambda = (3 + \sqrt{5})/2 > 1$, and 
$E^s(x)$ is the contracting eigenspace for $A$ corresponding to the smaller eigenvalue $\lambda^{-1} = (3 - \sqrt{5})/2 < 1$.  In this example, we can choose $n=1$ to verify that  uniform hyperbolicity holds in the definition.  The Lyapunov exponents of this cocycle are
$\pm\log\lambda$.    

The map $g_\epsilon\colon \TT^2\to \TT^2$ given by
\begin{equation}\label{e=gepsilon} g_\epsilon(x,y) = \left(2x + y + \epsilon \sin(2\pi (x + y)),  x+y  \right)
\end{equation}
is a $C^1$ small perturbation of $f_A$ if $\epsilon$ is sufficiently small, and so $g_\epsilon$ is Anosov for small $\epsilon$.   Moreover, since $\det D_{(x,y)} g_\epsilon \equiv 1$, the
map $g_\epsilon$ preserves the area $\mu$ on $\TT^2$; we shall see in the next section that $g_\epsilon$ is
ergodic with respect to $\mu$.   The  two Lyapunov exponents of $g_\epsilon$ with respect to the ergodic measure $\mu$ are  $\pm\log\lambda_\epsilon$, where $\lambda_\epsilon <\lambda$ for
$\epsilon\neq 0$ small. There are several ways to see this: one way to prove it is to compute directly using the ergodic theorem (Theorem~\ref{t=ergodic}) that $\epsilon\mapsto\lambda_\epsilon$ is a smooth map whose local maximum is achieved at $\epsilon=0$. 

In contrast to $f$, the exponents  of $g_\epsilon$ are not constant on $\TT^2$ but depend on the invariant measure. 
For example, the averaged Dirac measures $\nu_{x_1} = \delta_{(0,0)}$ and
$\nu_{x_2} = \frac{1}{3} \left(\delta_{(\frac12, \frac12)} + \delta_{(0,\frac12)} + \delta_{(\frac12,0)}\right)$
corresponding to the fixed point $x_1 = (0,0)$ and the periodic point $x_2 = (\frac12,\frac12)$ respectively, are both invariant
and ergodic under $g_\epsilon$, for any $\epsilon$.  Direct computation with the eigenvalues of the matrices
$D_{x_1} g_\epsilon$ and  $D_{x_2} g_\epsilon^3$ shows that the Lyapunov exponents with respect to $\nu_{x_1}$ and $\nu_{x_2}$ are different for $\epsilon\neq 0$.

\subsection{Measurably (nonuniformly) hyperbolic cocycles}

We say that a cocycle $A$ over $(\Omega,\mu,f)$ is {\em measurably hyperbolic} if for $\mu$-a.e. point $\omega\in \Omega$, the exponents $\chi_j(\omega)$ are all nonzero.  Since the role played by the measure is important in this definition, we sometimes say that $\mu$ is a hyperbolic measure for the cocycle $A$, or $A$ is hyperbolic with respect to $\mu$.

Uniformly hyperbolic cocycles over a homeomorphism $f$ are hyperbolic with respect to {\em any} $f$-invariant
probability measure $\mu$ (exercise).  An equivalent definition of measurable hyperbolicity that neatly parallels the uniformly hyperbolic condition is: $A$ is hyperbolic with respect to the $f$-invariant measure $\mu$ if there exist a set $X\subset \Omega$ with $\mu(X) = 1$ and  splittings $\RR^d = E^u(\omega)\oplus E^s(\omega)$ depending measurably on $\omega\in X$,  such that for every $\omega\in X$  there exists an integer $n = n(\omega)\geq 1$ such that conditions (i)-(iii) in Definition~\ref{d=unifhyp} hold.   The splitting into unstable and stable spaces is not necessarily continuous (or even globally defined on $\Omega$) and the amount of time $n(\omega)$ to wait for doubling in length to occur depends on $\omega$; for these reasons, measurably hyperbolic cocycles are often referred to as ``nonuniformly hyperbolic."
\footnote{The terminology is not consistent across fields.  In smooth dynamics, a cocycle over a measurable system that is measurably hyperbolic is called nonuniformly hyperbolic, whether it is uniformly hyperbolic or not.  In the spectral theory community, a  cocycle is called nonuniformly hyperbolic if it is measurably hyperbolic but {\em not} uniformly hyperbolic.}   
In this nonuniform setting it is possible for a cocycle to be hyperbolic with respect one invariant measure, but not another.

\begin{figure}[ht]
\begin{center}
\includegraphics[scale=.25 ]{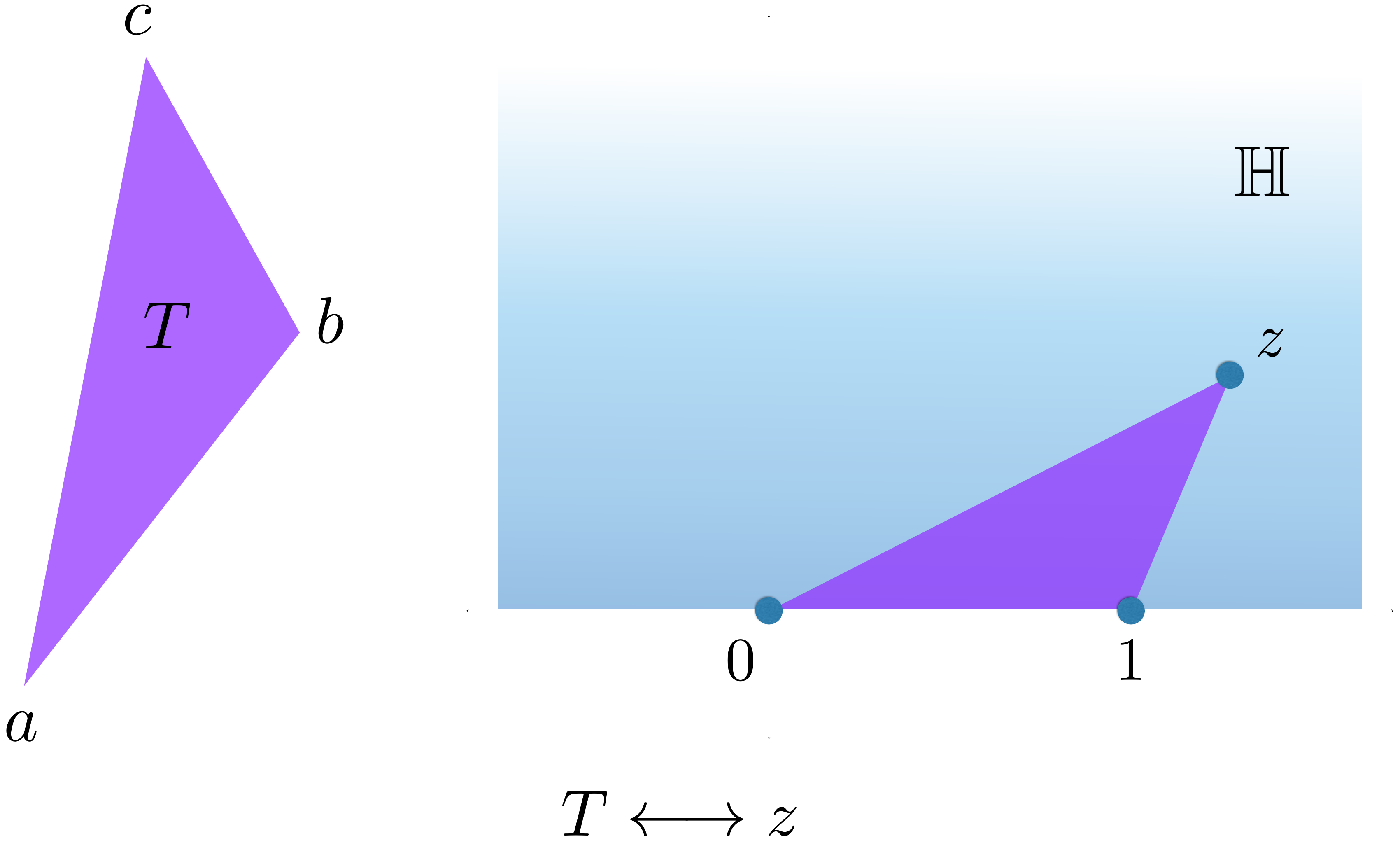}
\end{center}
\caption{The upper half plane $\HH$ is the space of marked triangles, up to Euclidean similarity. The triangle $T$ on the left corresponds to the point $z\in \HH$ on the right.}
\end{figure}

A measurably hyperbolic cocycle lurks behind random barycentric subdivision.
The random process generating the triangles in iterated barycentric subdivision can be encoded in a cocycle as follows. First, we identify upper half plane $\HH\subset \CC$  with the space of marked triangles (modulo Euclidean similarity) by sending a triangle $T$ with vertices cyclically labeled $a,b,c$ to a point $z\in \HH$ by rescaling, rotating and translating, sending $a$ to $0$, $b$ to $1$ and $c$ to $z$. See Figure 4.

\begin{figure}[ht]
\begin{center}
\includegraphics[scale=.27 ]{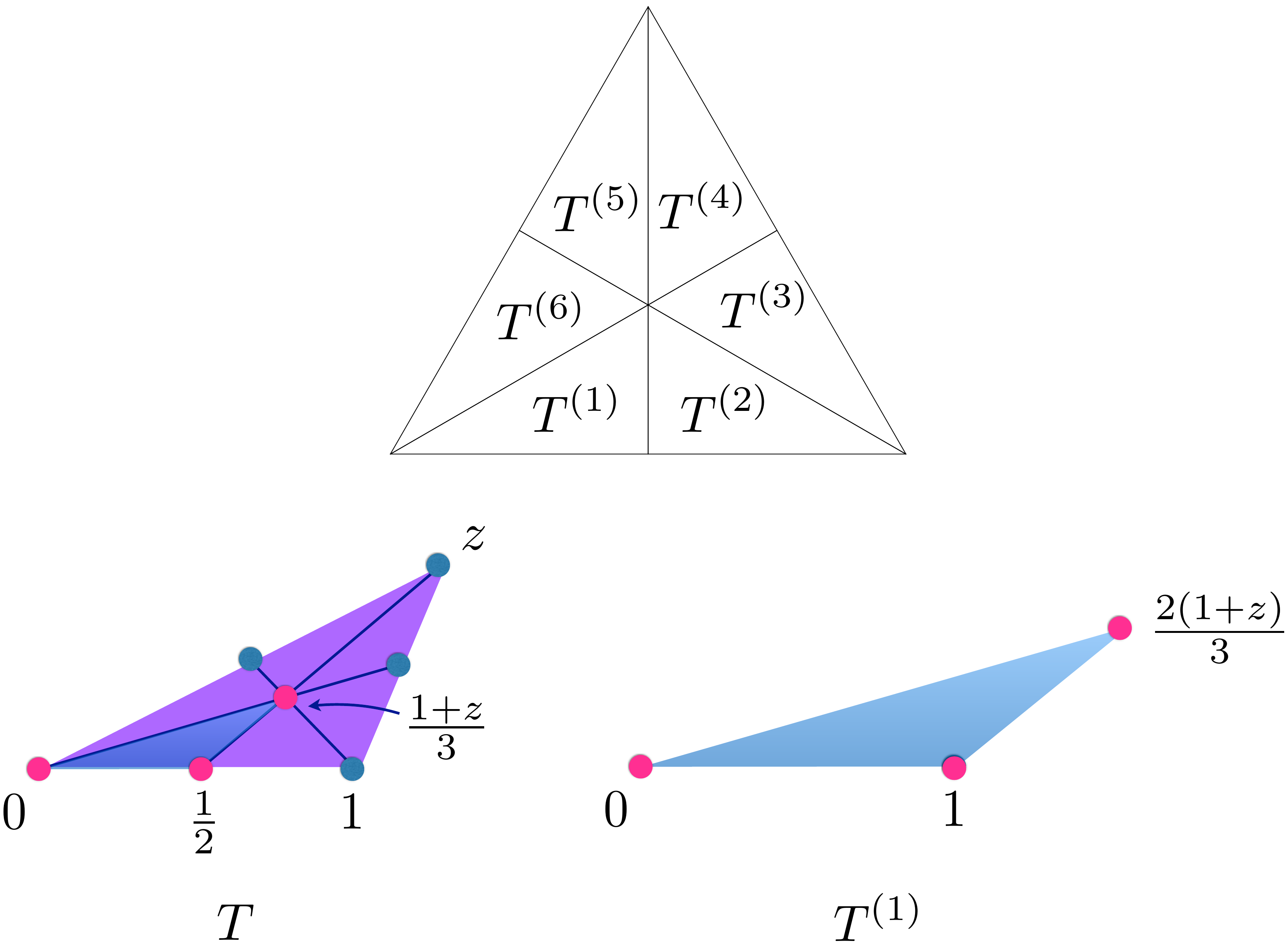}
\end{center}
\caption{A M\"obius transformation that selects the first triangle in barycentric subdivision.}
\end{figure}

Labelling cyclically the triangles $T^{(1)}, T^{(2)},\ldots, T^{(6)}$ in the subdivision as in Figure 5, the M\"obius transformation $B(z) = 2(z+1)/3$ sends the marked triangle $T$ to the marked triangle $T^{(1)}$.
The involutions  pictured in Figure 6  generate the symmetric group $S_3$, whose nontrivial elements we label $P_1, \ldots, P_5$. 
A bit of thought shows that the transformations $B, BP_1, \ldots , BP_5$ are the 6 maps of the plane sending $T$ to the  rescaled triangles in the subdivision.

Fixing an identification of the lower half plane with the upper half plane via $z\mapsto \bar z$, the action of these $6$ transformations on $\HH$ are identified with the projective action of $6$ elements  $A_1,\ldots,A_6$ of $PGL(2,\RR)$, where $B$ is identified with 
$\begin{pmatrix} 2/\sqrt{6} &2/\sqrt{6}\\0&3/\sqrt{6} \end{pmatrix}$, $P_1$ is identified with  $\begin{pmatrix} 1&0\\1&-1 \end{pmatrix}$, $P_2$ is identified with  $\begin{pmatrix} 0&1\\1&0 \end{pmatrix}$, and $P_3$ is identified with  $\begin{pmatrix} -1&1\\0&1 \end{pmatrix}$.
\begin{figure}[ht]
\begin{center}
\includegraphics[scale=.27 ]{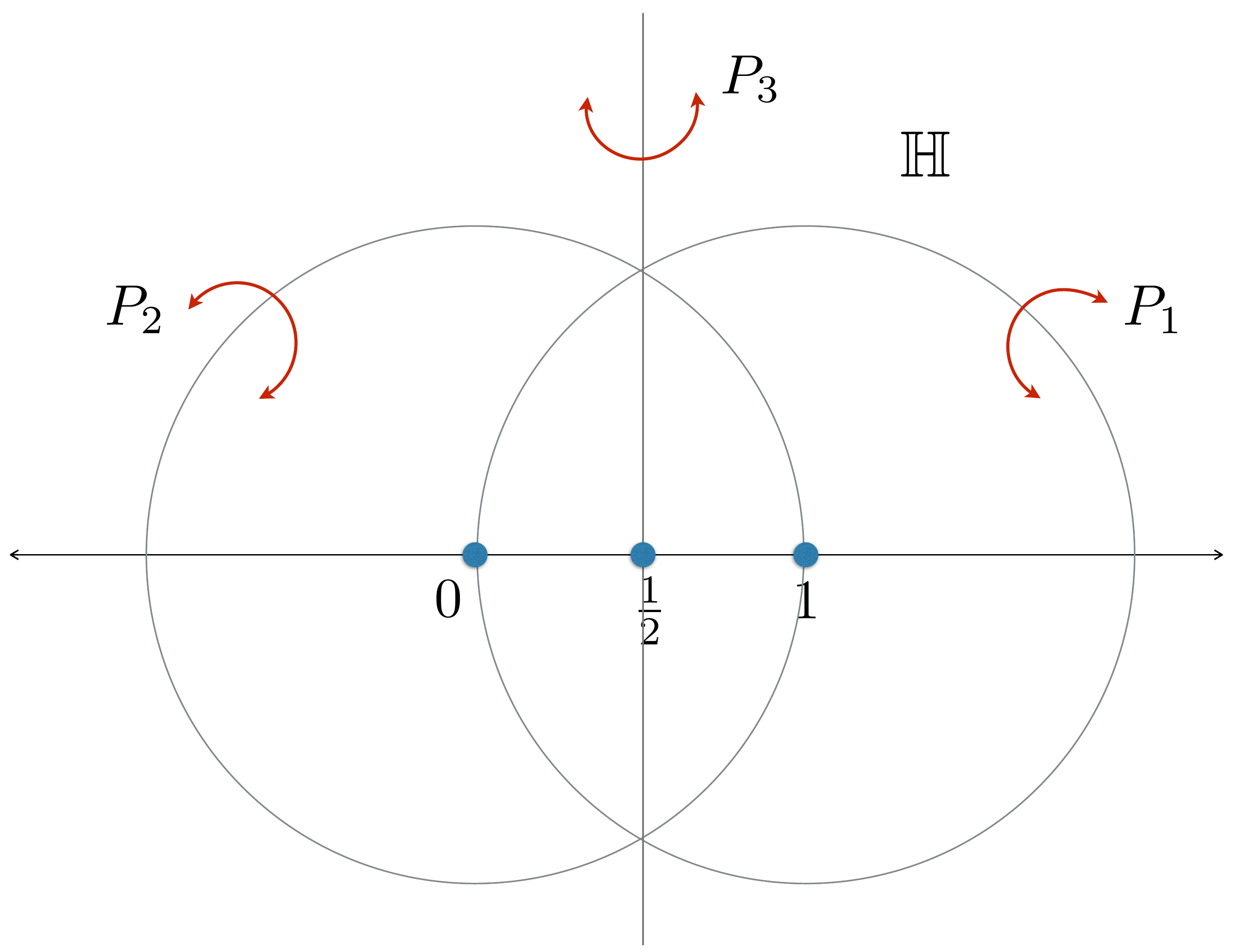}
\end{center}
\caption{Three involutions generating the symmetric group $S_3$: inversion in the circles $|z| = 1$ and $|z-1|=1$, and reflection across the line $\hbox{Re}(z) = \frac12$.}
\end{figure}

Thus random barycentric subdivision is governed by a random matrix cocycle over a Bernoulli shift.  If repeated rolls of the die generate the sequence  $\omega_1, \omega_2,\ldots$ with $\omega_i\in \{1,\ldots, 6\}$, then 
the $n$th triangle  $T_n(\omega)$ generated is the projective image of $\begin{pmatrix}i\\1\end{pmatrix}$ under $A^{(n)}(\omega)$, and an exercise shows
that the aspect ratio of $T_n(\omega)$ is given by:
 \[\alpha(T_n(\omega)) = \|(0,1)\cdot A^{(n)}(\omega)\|^{-2}.
\]
We thus obtain the formula
 \[\lim_{n\to\infty} \frac{1}{n}\log\alpha(T_n(\omega)) =  -2 \lim_{n\to\infty} \frac{1}{n}\log \left\| (0,1)\cdot A^{(n)}(\omega)\right\|,\]
and out pops the Lyapunov exponent $\chi = \lim_{n\to\infty} \frac{1}{n}\log \left\| (0,1)\cdot A^{(n)}(\omega)\right\|$.\footnote{As with eigenvalues, the Lyapunov exponents for left and right matrix multiplication coincide.}

Numerical simulation gives the value $\chi \approx  0.0446945$, but the fact that this number is positive follows from a  foundational result of Furstenberg (stated precisely as a special case in Theorem~\ref{t=furstenberg}) that  underlies some of Avila's results as well.  The upshot is that a random product of matrices
in $SL(2,\RR)$  (or $PGL(2,\RR)$)  cannot have exponents equal to 0, except by design.  In particular,
 if the matrices do not simultaneously preserve a collection of one or two lines, and the group generated by the matrices is not compact, then the exponents with respect to any nontrivial Bernoulli measure will be nonzero. These conditions are easily verified for the barycentric cocycle.   Details of this argument about barycentric subdivision can be found in \cite{McMo} and the related paper \cite{Barany:Beardon:Carne}.

The barycentric cocycle is {\em not} uniformly hyperbolic, as can be seen  by examining the sequence of triangles generated by subdivision on Figure 2: at each stage it is always possible to choose some triangle with aspect ratio bounded below, even though most triangles in a subdivision will have smaller aspect ratio than the starting triangle.
For products of matrices in $SL(d,\RR)$, uniform hyperbolicity must be carefully engineered, but for random products, measurable hyperbolicity almost goes without saying.  

A longstanding problem in smooth dynamics is to understand which diffeomorphisms have  hyperbolic derivative cocycle with respect to some natural invariant measure, such as volume (See \cite{BV1}).  Motivating this problem is the fact  that measurable hyperbolicity produces interesting dynamics, as we explain in the next section.

\section{Smooth ergodic theory }\label{s=lyaperg}

Smooth ergodic theory studies the dynamical properties of smooth maps from a statistical point of view.   A natural object of study is a measure-preserving system $(M, \hbox{vol}, f)$, where $M$ is a smooth, compact manifold without boundary equipped with a Riemannian metric, $\hbox{vol}$ is the volume measure of this metric, normalized so that $\hbox{vol}(M)=1$,  and $f\colon M\to M$ is a diffeomorphism preserving $\hbox{vol}$.  It was in this context that Boltzmann originally hypothesized ergodicity for ideal gases in the 1870's.  Boltzmann's non-rigorous formulation of ergodicity was close in spirit to the following statement of the pointwise ergodic theorem for diffeomorphisms: 

\begin{theorem}\label{t=ergodic} If $f$ is ergodic with respect to volume, then its orbits are equidistributed, in the following sense: for almost every $x\in M$, and any continuous function $\phi\colon M\to\RR$:
\begin{equation}\label{e=birklimit} \lim_{n\to\infty}\frac{1}{n}\left(\phi(x) + \phi(f(x))  + \cdots \phi(f^{n-1}(x)) \right) = \int_M \phi \, d\hbox{vol}.
\end{equation}
\end{theorem}

As remarked previously, an example of an ergodic diffeomorphism is the rotation $f_\alpha$ on $\RR/\ZZ$, for $\alpha$ irrational.   In fact this transformation has a stronger property of unique ergodicity, which is equivalent to the property that  the limit in (\ref{e=birklimit}) exists for {\em every} $x\in \RR/\ZZ$.\footnote{This is a consequence of Weyl's equidistribution theorem and can be proved using elementary analysis. 
See, e.g. \cite{Helson}. }  While unique ergodicity is a strong property, the ergodicity of irrational rotations is fragile; the ergodic map $f_\alpha$ can be perturbed to obtain the non-ergodic map $f_{p/q}$, where $\alpha\approx p/q$.

Another example of an ergodic diffeomorphism,  at the opposite extreme of the rotations in more than one sense, is the automorphism  $f_A$ of the 2-torus induced by multiplication by the matrix $A=\begin{pmatrix} 2&1\\1&1\end{pmatrix}$ with respect to the area $\mu$.    In spirit, this example is closely related to the Bernoulli shift, and in fact its orbits can be coded in such a way as to
produce a measure-preserving isomorphism with a Bernoulli shift. As observed in the previous section, ergodicity of this map 
can be proved using Fourier analysis, but there is a much more robust proof, due to Anosov \cite{Anosov},
who showed that any $C^2$ Anosov diffeomorphism preserving volume is ergodic with respect to volume.

\subsection{Ergodicity of Anosov diffeomorphisms and Pesin theory}

Anosov's proof of ergodicity has both analytic and geometric aspects.  For the map $f_A$, it follows several steps:
\begin{enumerate} 
\item The expanding and contracting subbundles $E^u$ and $E^s$ of the splitting $T\TT^2 = E^u\oplus E^s$ are tangent to foliations
$\cW^u$ and $\cW^s$ of $\TT^2$ by immersed lines.  These lines are parallel to the expanding and contracting eigendirections of $A$ and wind densely around the torus, since they have irrational slope.  The leaves of this pair of foliations are perpendicular to each other, since $A$ is symmetric.
\item  A clever application of the pointwise ergodic theorem (presented here as a special case in Theorem~\ref{t=ergodic}) shows that any $\phi\in L^2(\TT^2,\mu)$ satisfying $\phi\circ f = \phi$ is,
up to a set of area 0, constant along leaves of the foliation $\cW^u$,  and (again, up to a set of area 0) constant along 
leaves of $\cW^s$. This part of the argument goes back to Eberhard Hopf's study of geodesics in negatively curved surfaces in the 1930's.
\item Locally, the pair of foliations $\cW^u$ and $\cW^s$ are just (rotated versions of) intervals parallel to the $x$ and $y$ axes. In these rotated coordinates, $\phi(x,y)$ is a measurable function constant a.e.  in $x$ and constant a.e. in $y$.  Fubini's theorem then implies that such a $\phi$ must be constant a.e.  This conclusion holds in local charts, but since $\TT^2$ is connected, $\phi$ must be constant.
\item Since any $f$-invariant function $\phi\in L^2(\TT^2,\mu)$ is constant almost everywhere with respect to $\mu$, we conclude that $f_A$ is ergodic with respect to $\mu$.
\end{enumerate}

The same proof works for any smooth, volume-preserving Anosov diffeomorphism --- in particular,  for  the maps $g_\epsilon$ defined in (\ref{e=gepsilon}) --- if one modifies the steps appropriately.  The foliations by parallel lines in step (1)  are replaced by foliations by smooth curves (or submanifolds diffeomorphic to some $\RR^k$, in higher dimension).  Step (2) is almost the same, since it uses only volume preservation and the fact that the leaves of $\cW^u$ and $\cW^s$ are expanded and contracted, respectively. Step (3) is the most delicate to adapt and was Anosov's great accomplishment, and Step (4) is of course the same.
 
The ergodicity of Anosov diffeomorphisms is also a consequence of the much stronger property of being measurably encoded by a Bernoulli shift.
This so-called {\em Bernoulli property} implies that Anosov diffeomorphisms are mixing with respect to volume, meaning that for any $L^2$ function $\phi$, we have
$\int_M \phi \cdot \phi\circ f^n \to (\int_M \phi)^2$ as $n\to \infty$.  Visually, sets are mixed up by Anosov diffeomorphisms:  see Figure 4, which explains why $f_A$ is sometimes referred to as the ``cat map."
\begin{figure}[h]
\begin{center}
\includegraphics[scale=0.5]{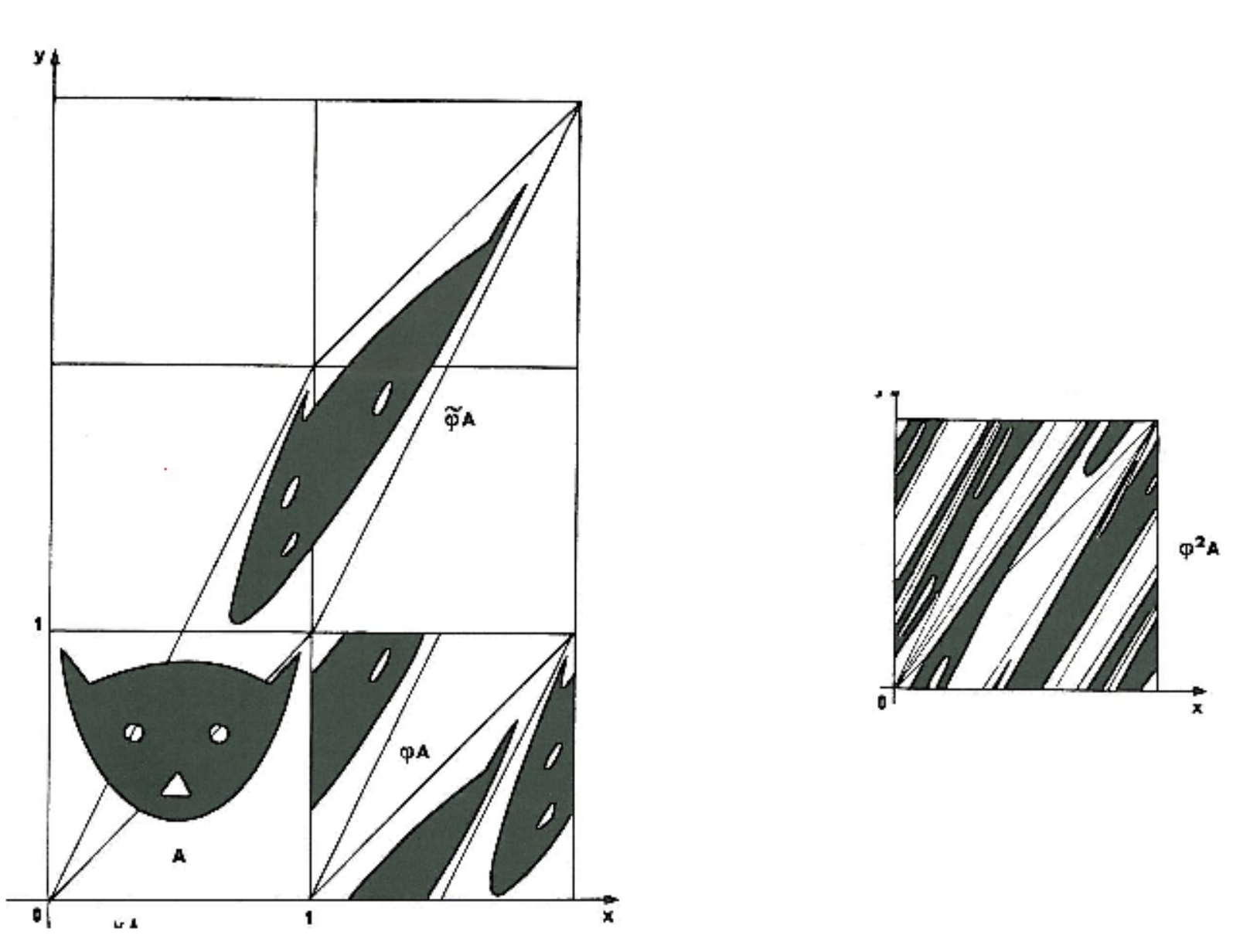}
\end{center}
\caption{The action of $f_A$ on a  cat, from \cite{AA}. A cat is drawn in a square fundamental domain for $\TT^2$ at the lower left.  Its image under $A$ is shown in the parallelogram, and it is reassembled  into another fundamental domain to show its image under $f_A$.  The image of the cat under $f_A^2$ is depicted at right. }
\end{figure}
The map $g_\epsilon$ will similarly do a number on a cat.   The proof of the Bernoulli property for Anosov diffeomorphisms builds the Anosov-Hopf proof of ergodicity.  The Anosov-Hopf argument has the additional advantage that it can be adapted to prove ergodicity for systems that are not Bernoulli.

As explained in Subsection~\ref{ss=Anosov}, any $C^1$-small perturbation of an Anosov diffeomorphism is Anosov, and  $C^2$ Anosov diffeomorphisms preserving volume are ergodic.
Hence volume-preserving  $C^2$ Anosov diffeomorphisms  are  {\em stably ergodic:} the ergodicity cannot be destroyed by a $C^1$ small perturbation, in marked contrast with the irrational rotation $f_\alpha$. 
 
The Anosov condition thus has powerful consequences in smooth ergodic theory.  But, like uniformly hyperbolic matrix products, Anosov diffeomorphisms are necessarily contrived.  In dimension 2, the only surface supporting an Anosov  diffeomorphism is the torus $\TT^2$, and 
conjecturally, the only manifolds supporting Anosov diffeomorphisms belong to a special class called the infra-nilmanifolds.  
On the other hand, {\em every} smooth manifold supports a volume-preserving diffeomorphism that is hyperbolic with respect to volume, as was shown by Dolgopyat and Pesin in 2002 \cite{DP}.

In the 1970's Pesin \cite{Pe} introduced a significant innovation in smooth ergodic theory: a nonuniform, measurable analogue of the Anosov-Hopf theory.  Under the assumption that a volume preserving diffeomorphism $f$ is hyperbolic with respect to volume, Pesin showed that volume has at most countably many ergodic components with respect to $f$.  Starting with Oseledets's theorem, and repeatedly employing Lusin's theorem that every measurable function is continuous off of a set of small measure, Pesin developed an ergodic theory of smooth systems that has had numerous applications.  Some limitiations of Pesin theory are: first, that it begins with the hypothesis of measurable hyperbolicity, which is a condition that is very hard to verify except in special cases, and second, without additional  input, measurable hyperbolicity does not imply ergodicity, as the situation of infinitely many ergodic
components can and does occur \cite{DHP}.

\subsection{Ergodicity of ``typical" diffeomorphisms}

The question of whether ergodicity is a common property among volume-preserving diffeomorphisms of a compact manifold $M$ is an old one,
going back to Boltzmann's ergodic hypothesis of the late  19th Century.
We can formalize the question by fixing a differentiability class $r\in[1,\infty]$ and considering the set $\Diff^r_{\hbox{\tiny vol}}(M)$ of 
$C^r$, volume-preserving diffeomorphisms of $M$.  This is a topological space in the $C^r$ topology, and we say that a property holds {\em generically in $\Diff^r_{\hbox{\tiny vol}}(M)$} (or {\em $C^r$ generically}, for short)
if it holds for all $f$ in a countable intersection of open and dense subsets of $\Diff^r_{\hbox{\tiny vol}}(M)$.\footnote{Since $\Diff^r_{\hbox{\tiny vol}}(M)$ is a Baire space, properties that hold generically hold for a dense set, and two properties that hold generically separately hold together generically.}

\begin{figure}[h]
\begin{center}
\includegraphics[scale=0.27]{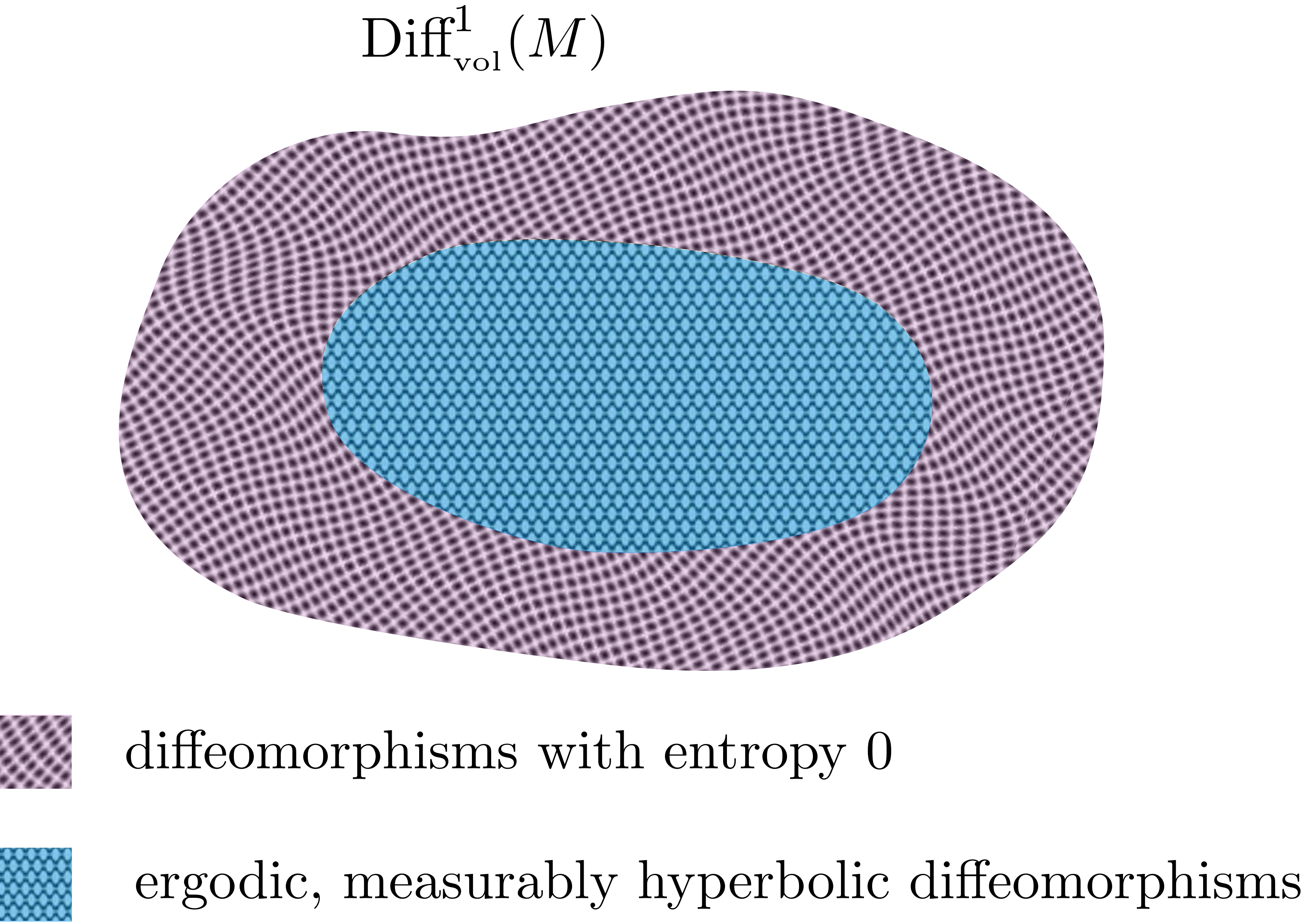}
\end{center}
\caption{Generically, positive entropy implies ergodicity and measurable hyperbolicity.}
\end{figure}

Oxtoby and Ulam \cite{OU} proved in 1939 that a generic volume-preserving  {\em homeomorphism} of a compact manifold is ergodic.  At the other extreme, KAM (Kolmogorov-Arnol'd-Moser) theory, introduced by Kolmogorov in the 1950's \cite{Ko},  implies that ergodicity is {\em not} a dense property, let alone a generic one, in $\Diff^\infty_{\hbox{\tiny vol}}(M)$, if $\hbox{dim}(M)\geq 2$.  The general question of whether ergodicity is generic in  $\Diff^r_{\hbox{\tiny vol}}(M)$ remains open for $r\in [1,\infty)$, but we now have a complete answer for any manifold when $r=1$
under the assumption of {\em positive entropy}. Entropy is a numerical invariant attached to a measure preserving system that measures the complexity of orbits.  The rotation $f_\alpha$ has entropy $0$; the Anosov map $f_A$ has positive entropy $\log(\lambda)$.  By a theorem of Ruelle, positivity of entropy means that there is {\em some} positive volume subset of $M$ on which the Lyapunov exponents are nonzero in {\em some} directions.

\begin{theorem}[Avila, Crovisier, Wilkinson \cite{ACW}]\label{t=acw} Generically in $\Diff^1_{\hbox{\tiny vol}}(M)$, positive entropy implies ergodicity and moreover measurable hyperbolicity with respect to volume.
\end{theorem}

See Figure 8. This result was proved in dimension 2 by Ma\~n\'e-Bochi \cite{M,B} and dimension $3$ by M.A. Rodriguez-Hertz \cite{jana}.
Positive entropy is an {\em a priori} weak form of chaotic behavior that can be confined to an invariant set of very small measure, with trivial dynamics on the rest of the manifold.  Measurable hyperbolicity, on the other hand, means that at almost every point {\em all} of the Lyapunov exponents of the derivative cocycle $Df$ are nonzero.
Conceptually, the proof divides into two parts:
\begin{enumerate}
\item  $C^1$ generically, positive entropy implies nonuniform hyperbolicity.  One needs to go from some nonzero exponents on some of the manifold to  all nonzero exponents on almost all of the manifold.  Since  the cocycle and the dynamics are intertwined, carrying this out is a delicate matter. This relies on the fact that the $C^1$ topology
is particularly well adapted to the problem.  On the one hand, constructing $C^1$-small perturbations with a desired property is generally much easier than constructing $C^2$ small perturbations with the same property.  On the other hand, many useful dynamical properties such as uniform hyperbolicity are $C^1$ open.
\item  $C^1$ generically, measurable hyperbolicity (with respect to volume) implies ergodicity.  This argument uses Pesin theory but adds some missing input needed to establish ergodicity.  This input holds $C^1$ generically.   For example, the positive entropy assumption generically implies existence of a {\em dominated splitting}; this means that generically, a positive entropy diffeomorphism  is something intermediate between an Anosov diffeomorphism and a nonuniformly hyperbolic one.  There is a {\em continuous} splitting $TM = E^u\oplus E^s$, invariant under the derivative, such that for almost every $x\in M$, there exists an $n= n(x)$ such that
for every $\xi^u\in E^u(x)$,  $\|D_xf^n(\xi^u)\| \geq 2\|\xi^u\|$, and for every $\xi^s\in  E^s(x)$, $\|D_xf^{-n}(\xi^s)\| \geq 2\|\xi^s\|$.
\end{enumerate}
The proof incorporates techniques from several earlier results, most of which have been proved in the past 20 years \cite{AB,Boc2, BV2, SW}.
Also playing an essential technical role in the argument is a regularization theorem of Avila: every $C^1$ diffeomorphism that preserves volume can be $C^1$ approximated by a $C^2$ volume-preserving diffeomorphism \cite{A}.  The fact that this regularization theorem was not proved until recently highlights the difficulty in perturbing the derivative cocycle to have desired properties: you can't change $Df$ without changing $f$ too (and {\em vice versa}). This is why completely general results analogous to Furstenberg's theorem for random matrix products are few and far between for diffeomorphism cocycles.

\section{Translation surfaces}

A flat surface is any closed surface that can be obtained by gluing together finitely many parallelograms in $\RR^2$ along coherently oriented parallel edges, as in Figure 9.
\begin{figure}[h]
\begin{center}
\includegraphics[scale=0.28]{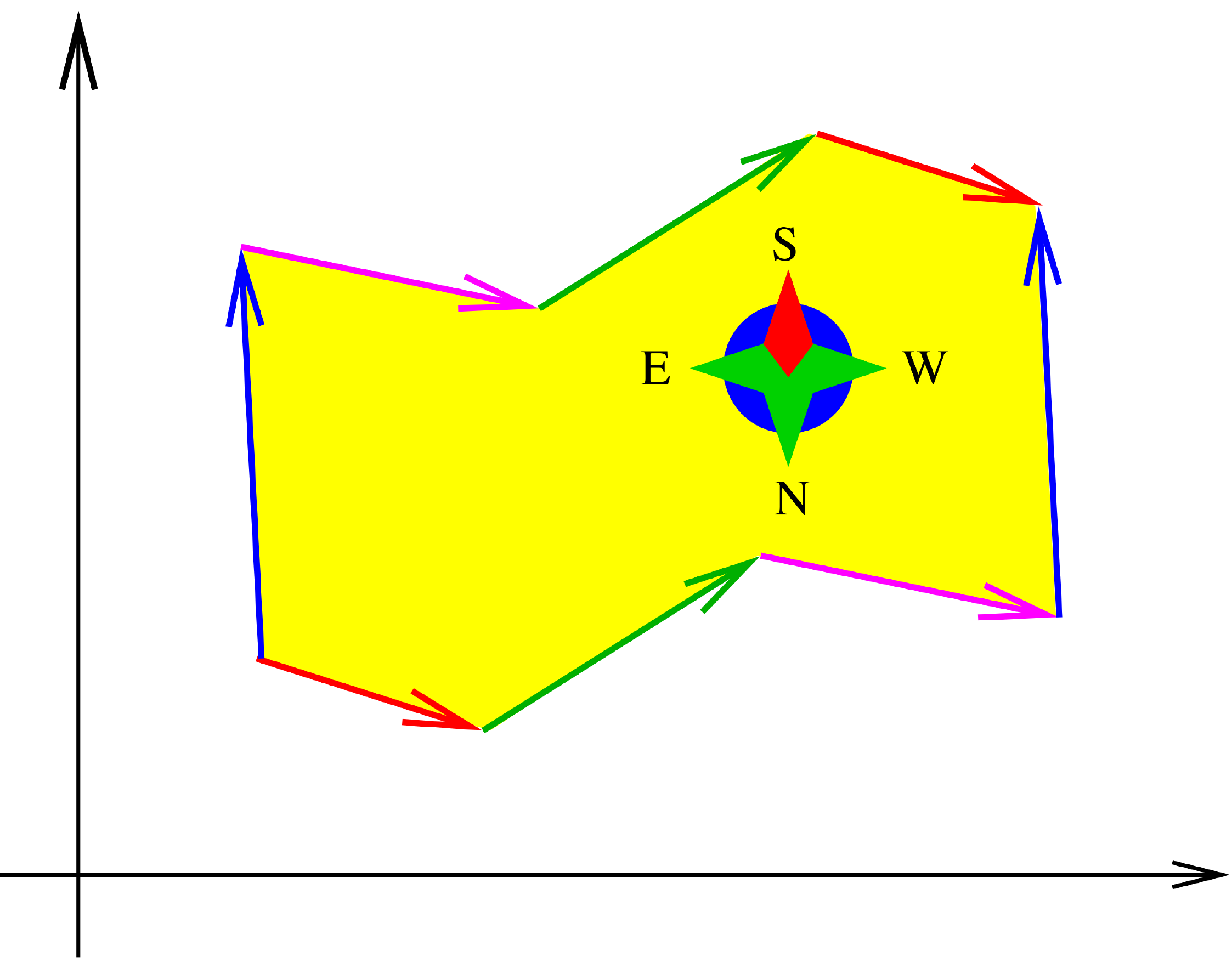}
\end{center}
\caption{A flat surface with a distinguished ``South," also known as a translation surface (courtesy Marcelo Viana).
Parallel edges of the same color are identified.}
\end{figure}
Two flat surfaces are equivalent if one can be obtained from the other by cutting,  translating,  and rotating.
A  {\em translation} surface is a flat surface that comes equipped with a well-defined, distinguished vertical,  ``North" direction (or, ``South" depending on your orientation).  Two  translation surfaces are equivalent if one 
can be obtained from the other by cutting and translating  (but {\em not} rotating).

Fix a translation
surface $\Sigma$ of genus $g>0$.
 If one picks an angle $\theta$ and a point $x$ on $\Sigma$,  and follows the corresponding straight ray through  $\Sigma$, there are two possibilities: either it terminates in a corner, or it can be continued for all time.  Clearly for any $\theta$, and almost every starting point (with respect to area), the ray will continue forever.  If it continues forever, either it returns to the initial point and direction and produces a closed curve, or it continues on a parallel course without returning.  A version of the Pigeonhole Principle for area  (Poincar\'e recurrence) implies that for almost every point and starting direction, the line will come back arbitrarily close to the starting point.

\begin{figure}[h]
\begin{center}
\includegraphics[scale=0.28]{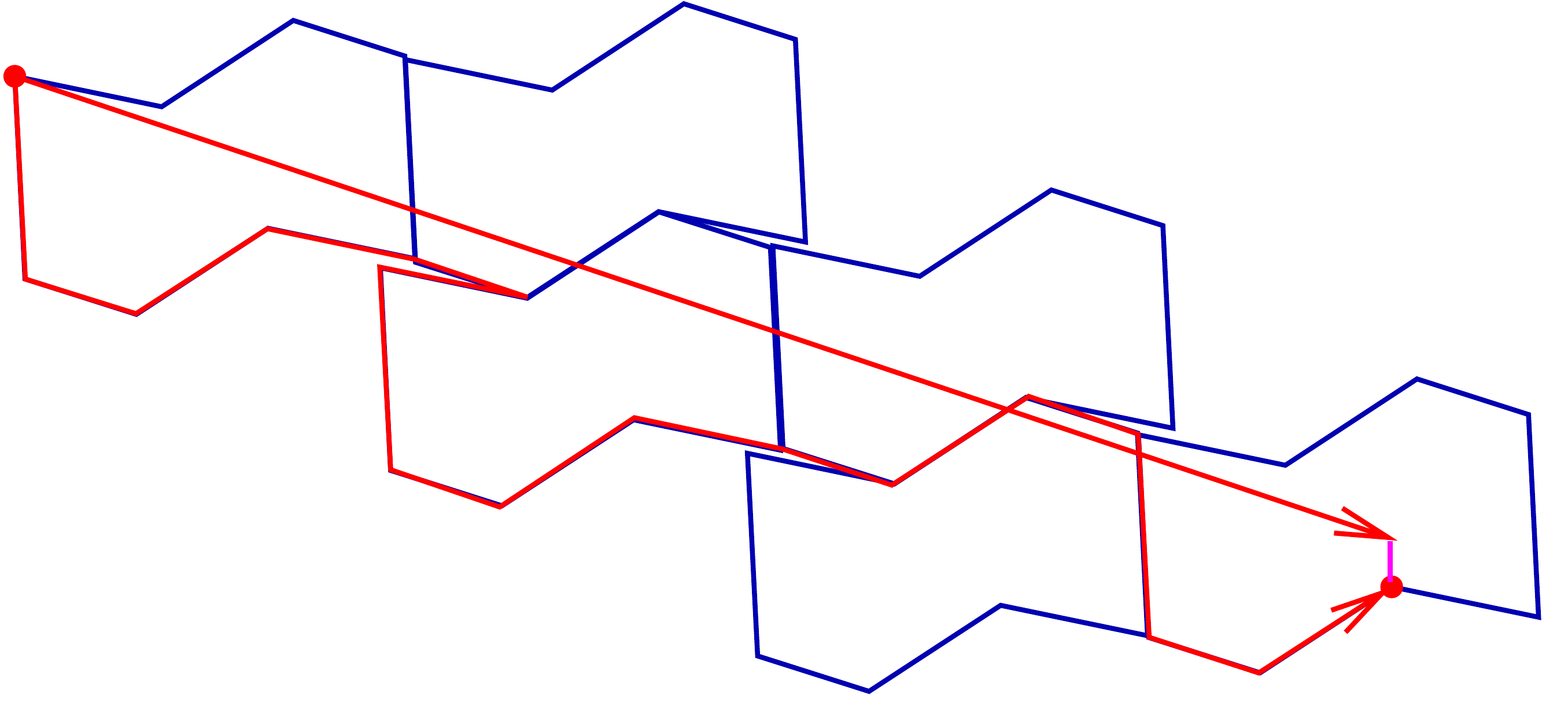}
\end{center}
\caption{Closing up a ray  that comes back close to itself (courtesy Marcelo Viana)}
\end{figure}

Kerckhoff-Masur-Smillie \cite{kms} proved more: for a fixed $\Sigma$, and almost every $\theta$, the ray through any point $x$ is dense in $\Sigma$, and in fact is equidistributed with respect to area.
Such a direction $\theta$ is called {\em uniquely ergodic}, as it is uniquely ergodic in the same sense that $f_\alpha$ is, for irrational $\alpha$.
Suppose we start with a uniquely ergodic direction and wait for the successive times that this 
ray returns closer and closer to itself.  This produces a sequence of closed curves $\gamma_n$
which produces a sequence of cycles $[\gamma_n]$ in homology $H_1(\Sigma, \ZZ) \simeq \ZZ^{2g}$.

Unique ergodicity of the direction $\theta$ implies that there is a unique $c_1\in H_1(\Sigma,\RR)$ such that for any starting point $x$:
$$
\lim_{n\to\infty} \frac{[\gamma_n]}{\ell(\gamma_n)} = c_1,
$$
where $\ell(\gamma)$ denotes the length in $\Sigma$ of the curve $\gamma$.

\begin{theorem} [Forni, Avila-Viana, Zorich \cite{Fo, AV, Zorich1, Zorich2}]\label{t=favc} Fix a topological surface $S$ of genus $g\geq 1$, 
and let $\Sigma$ be almost any translation surface modelled on $S$.\footnote{``Almost any" means with respect to the Lebesgue measure on possible choices of lengths and directions for the sides of the pentagon. This statement can be made more precise in terms of Lebesgue measure restricted to various strata in the moduli space of translation surfaces.} Then there exist real numbers $1 >\nu_2>\ldots>\nu_g>0$ and  a sequence of of subspaces
$L_1\subset L_2 \subset \cdots L_g$ of  $H_1(\Sigma,\RR)$  with $\dim(L_k) = k$ such that
for almost every $\theta$, for every $x$, and every $\gamma$ in direction $\theta$,
 the distance from $[\gamma]$ to $L_g$ is bounded, and
$$\limsup_{\ell(\gamma)\to \infty} \frac{\log \dist([\gamma], L_i)}{\log(\ell(\gamma))} = \nu_{i+1},$$
for all $i<g$. 

\end{theorem}
This theorem gives precise information about the way  the direction of $[\gamma_n]$ converges  to its asymptotic cycle $c_1$:  the convergence has a ``directional nature" much in the way
a vector $v\in \RR^d$ converges to infinity under repeated application of a matrix
\[
A =  \begin{pmatrix}  \lambda_1 & \ast & \cdots  & \ast \\
 0 & \lambda_2  & \cdots & \ast \\
0 & \cdots &\cdots &\ast\\
0 & 0 &\cdots & \lambda_d 
\end{pmatrix},
\]
with $\lambda_1 > \lambda_2 > \cdots \lambda_d > 1$.

\begin{figure}[h]
\begin{center}
\includegraphics[scale=0.28]{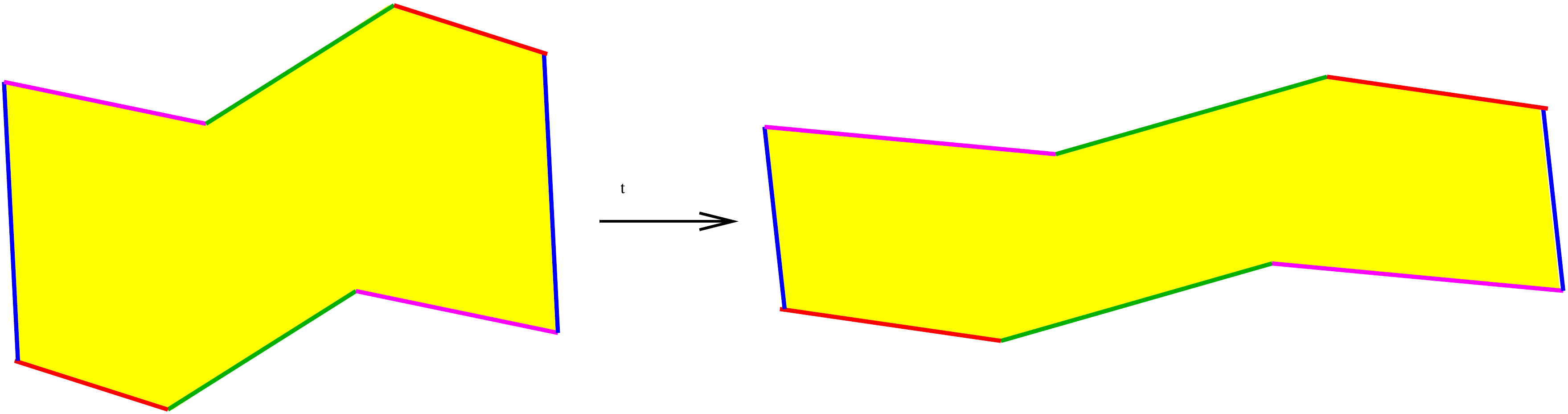}
\end{center}
\caption{A local picture of the Teichm\"uller flow (courtesy Marcelo Viana).}
\end{figure}

The numbers $\nu_i$ are the Lyapunov exponents of the {\em Kontsevich-Zorich} (KZ) cocycle over the 
{\em Teichm\"uller flow}.  The Teichm\"uller flow $\cF_t$  acts on the moduli space $\cM$ of translation surfaces 
(that is, translation surfaces modulo cutting and translation) by stretching in the East-West direction and contracting in the North-South direction. 
More precisely, if $\Sigma$ is a translation surface, then $\cF_t(\Sigma)$
is a new surface, obtained by transforming $\Sigma$ by the linear map $\begin{pmatrix}e^t & 0 \\ 0 & e^{-t} \end{pmatrix} $.
Since a stretched surface can often be reassembled to obtain a more compact one, it is plausible  that the Teichm\"uller flow has recurrent orbits (for example, periodic orbits). 
This is true  and  reflects  the fact that the flow $\cF_t$ preserves a natural volume that assigns finite measure to $\cM$.
The Kontsevich-Zorich
cocycle takes values in the symplectic group $Sp(2g,\RR)$ and captures homological data about
the cutting and translating equivalence on the surface.

 Veech proved that $\nu_2<1$ \cite{Veech},   Forni proved that $\nu_g > 0$ \cite{Fo},
and Avila-Viana proved that the numbers $\nu_2, \nu_3,\dots, \nu_{g-1}$ are all distinct \cite{AV}.  Zorich established the connection between exponents and the so-called deviation spectrum,
which holds in greater generality \cite{Zorich1, Zorich2}.  
\begin{figure}[ht]
\includegraphics[scale=.6]{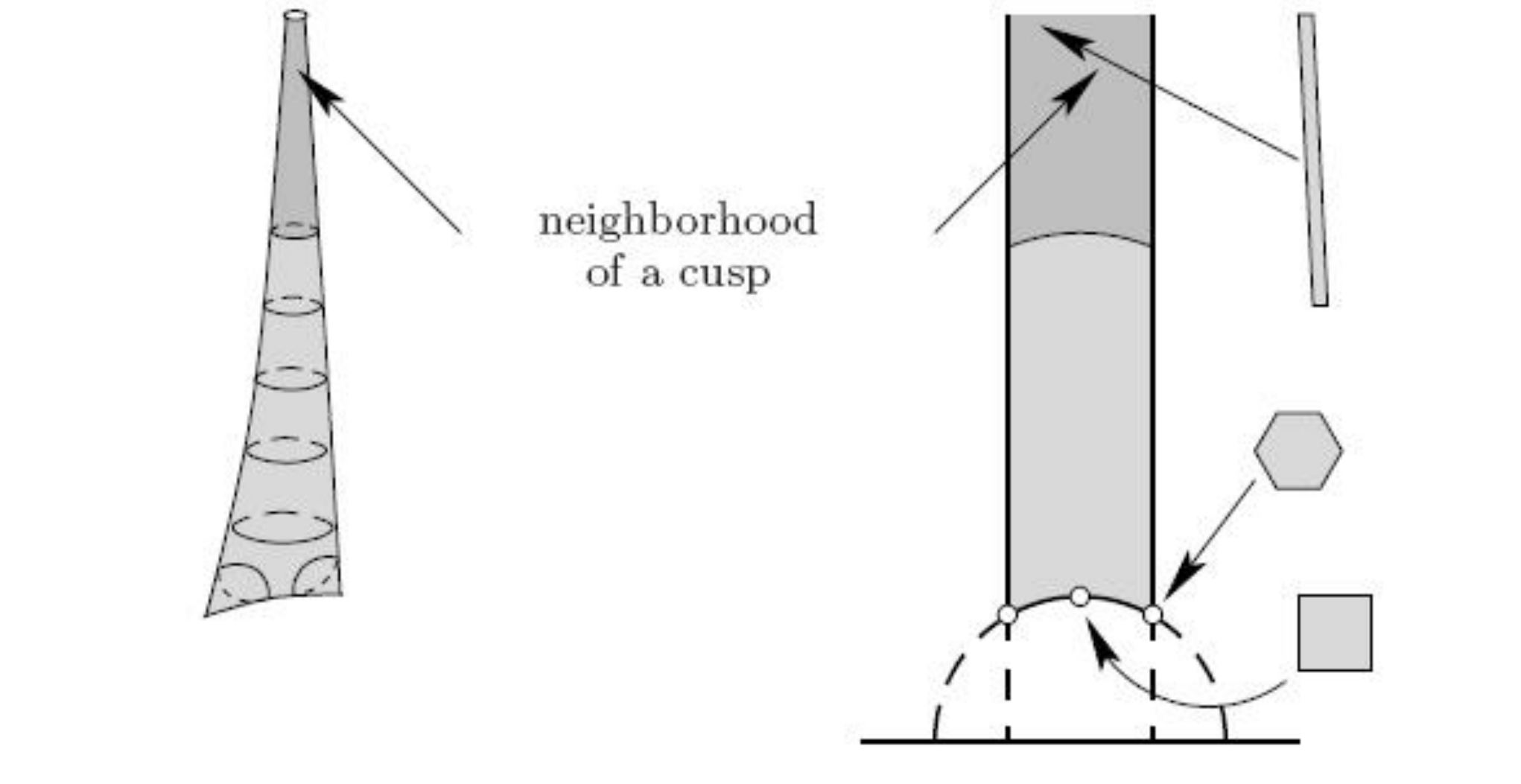}
\caption{The moduli space $\HH^2/PSL(2,\ZZ)$ of flat structures on the torus, which is a punctured sphere with two cone points, corresponding to the hexagonal and square lattices in $\RR^2$. Pictured at right is a fundamental domain in $\HH^2$.  (Courtesy Carlos Matheus.)}
\end{figure}
Many more things have been proved about the Lyapunov exponents of the KZ cocycle, and some of their values have been calculated which are (until recently, conjecturally) rational numbers!  See \cite{ekm, ce}.

Zorich's result reduces the proof  of Theorem~\ref{t=favc} to proving that the
the exponents  $\nu_1, \ldots, \nu_g$ are positive and distinct.
In the $g=1$ case where $\Sigma$ is a torus, this fact has a simple explanation. 
The moduli space $\cM$ is the set of all  flat structures on the torus (up to homothety), equipped with a direction. This is the quotient $SL(2,\RR)/SL(2,\ZZ)$, which is the unit tangent bundle of the modular surface $\HH/SL(2,\ZZ)$.  The  (continuous time) dynamical system  $\cF_t$ on  $\Omega$ is given by left multiplication by the matrix
$\begin{pmatrix}e^t & 0 \\ 0 & e^{-t} \end{pmatrix}$.  The cocycle is, in essence, the derivative 
cocycle for this flow (transverse to the direction of the flow)  This flow is {\em uniformly hyperbolic} (i.e. Anosov), and its exponents are  $\nu_1 = \log(e) = 1$ and $-\nu_1 = \log(e^{-1}) = -1$.  

The proof  for general translation surfaces that the exponents $\nu_1, \ldots, \nu_g$ are positive and distinct is considerably more involved.  We can nonetheless boil it down to some basic ideas. 

\begin{enumerate}
\item  The Teichm\"uller flow itself is nonuniformly hyperbolic with respect to a natural volume  (Veech \cite{Veech}), and can be coded in a way that the dynamics appear almost random.
\item Cocycles over perfectly random systems (for example i.i.d. sequences of matrices) have a tendency to have distinct, nonzero Lyapunov exponents.  This was first proved by Furstenberg in the $2\times 2$ case \cite{F} and later by  Guivarc'h-Raugi \cite{GR} (see also  \cite{GM}).
\item Cocycles over systems that are nonrandom, but sufficiently hyperbolic and with good coding, also tend to have distinct, nonzero Lyapunov exponents.  This follows from  series of results, beginning with Ledrappier in the $2\times 2$ case \cite{L}, and in increasing generality by Bonatti-Viana \cite{BoVi},  Viana \cite{Almost}, and Avila-Viana \cite{AVp}.
\end{enumerate}

\section{Hofstadter's butterfly}\label{s=Schrodinger}

Pictured in Figure 13 is the spectrum of the operator $H^{\alpha}_x\colon \ell^2(\ZZ)\to \ell^2(\ZZ)$ given by
\[ [H^{\alpha}_x u](n) = u(n+1) + u(n-1) + 2 \cos(2\pi (x + n\alpha)) u(n), \]
where $x$ is a fixed real number called the {\em phase}, and $\alpha\in [0,1]$ is a parameter called the {\em frequency}.  The vertical variable is  $\alpha$, and the horizontal variable is the spectral energy parameter $E$, which ranges in $[-4,4]$.  We can read off the spectrum of $H^{\alpha}_x$ by taking a horizontal slice at height $\alpha$; the black region is the spectrum.

\begin{figure}[h]
\begin{center}
\includegraphics[scale=0.75]{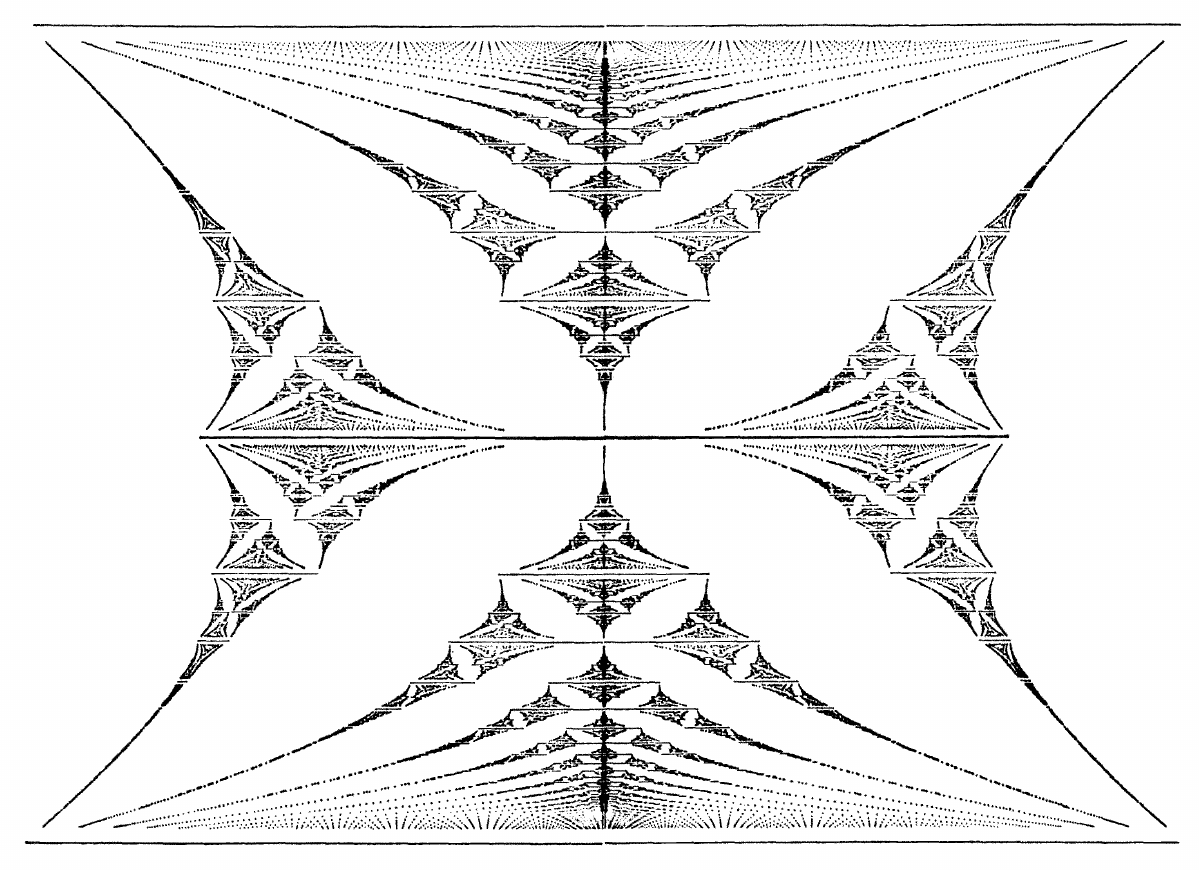}
\end{center}
\caption{Hofstadter's butterfly,  from \cite{Hof}.}
\end{figure}

In an influential 1976 paper, Douglas Hofstadter of {\em G\"odel, Escher, Bach} fame discovered this fractal picture while modelling the behavior of electrons in a crystal lattice under the force of a magnetic field \cite{Hof}. 
 This  operator plays a central role in the theory of the integer  quantum Hall effect developed by Thouless et al.,
 and, as predicted theoretically,  the butterfly has indeed appeared in von
Klitzing's QHE experiments. 
Similar operators are used in modeling 
graphene, and similar butterflies also appear in  graphene related
experiments (see, e.g. \cite{Detal}).

Some properties of the butterfly have been established rigorously.  For example, Avila and Krikorian proved:
\begin{theorem} [Avila-Krikorian, \cite{AK}] For every irrational
$\alpha\in [0,1]$, the $\alpha$-horizontal slice of the butterfly has measure $0$.
\end{theorem}
Their proof complements and thus extends the earlier result of Last  \cite{Last}, who proved the same statement, but for a full measure set of $\alpha$ satisfying an arithmetic condition.  In particular, we have:
\begin{corollary}\label{c=meas0} The  butterfly has measure $0$.
\end{corollary}
Other properties of the butterfly,  for example its Hausdorff dimension, remain unknown.

The connection between the spectrum of this operator and cocycles is an interesting one.
Recall the definition of the spectrum of $H_x^{\alpha}$:
$$
\sigma(H_x^\alpha):=\{E\in\CC: H_x^{\alpha}-E\cdot Id \mbox{ is not invertible}\}.
$$ 
The  eigenvalues are those $E$ so that the
eigenvalue equation $H^\alpha_{x}u=Eu$ admits $\ell^2(\ZZ)$ solutions, i.e. those $E$ such that
$H_x^{\alpha}-E\cdot Id$ is not injective.

The following simple observation is key. 
A sequence $(u_n:n\in\ZZ)\subset \CC^\ZZ$ (not necessarily in $\ell^2(\ZZ)$) solves $H_{x}^\alpha u=Eu$ if and only if
$$
A_E(f_\alpha ^{n}(x))\binom{u_{n}}{u_{n-1}}=\binom{u_{n+1}}{u_{n}},\ n\in\ZZ,
$$
where $f_\alpha\colon \RR/\ZZ\to \RR/\ZZ$ is the translation mentioned above, and
\begin{equation}\label{schrodinger-cocycle-map}
A_E(x)=\begin{pmatrix}E- 2\cos (2\pi x)& -1\\1& 0\end{pmatrix},
\end{equation}
which defines an $SL(2,\RR)$-cocycle over the rotation $f_\alpha$, an example of a {\em Schr\"odinger cocycle}.
Using the cocycle notation, we have
$$
A^{(n)}_E(x)\binom{u_{0}}{u_{-1}}=\binom{u_{n}}{u_{n-1}},\ n\in\ZZ.
$$

Now let's connect the properties of this cocycle with the spectrum of $H_{x}^\alpha$.
Suppose for a moment that the cocycle $A_E$ over $f_\alpha$ is uniformly hyperbolic, for some value of $E$. Then for every $x\in \RR/\ZZ$ there is a splitting $\RR^2 = E^u(x) \oplus E^s(x)$
invariant under cocycle, with vectors in $E^u(x)$ expanded under $A_E^{(mn)}(x)$, and vectors in
 $E^s(x)$ expanded under $A_E^{(-mn)}(x)$, both  by a factor of $2^m$.   Thus no solution $u$ to $H^\alpha_{x}u=Eu$  can be polynomially bounded simultaneously
in both directions, which implies $E$ is not an $\ell^2$ eigenvalue of $H^\alpha_{x}$.  It turns out
that the converse is also true, and moreover:
\begin{theorem}[R. Johnson, \cite{johnson}] \label{t.johnson}
If $\alpha$ is irrational, then for every $x\in [0,1]$:
\begin{equation}
\sigma(H_x^\alpha) =\{E: A_E \hbox{ is {\em not} uniformly hyperbolic over } f_\alpha\}.
\end{equation}
\end{theorem}
For irrational $\alpha$, we denote by $\Sigma_{\alpha}$ the spectrum of $\sigma(H_x^\alpha)$, which by Theorem~\ref{t.johnson} does not depend on $x$.  Thus for irrational $\alpha$, the set  $\Sigma_{\alpha}$ 
is the $\alpha$-horizontal slice of the butterfly.

 The butterfly is therefore both a dynamical picture and a spectral one.  On the one hand, it depicts the spectrum of a family of operators $H^\alpha_x$ parametrized by $\alpha$, and 
on the other hand, it depicts, within a $2$-parameter family of cocycles $\{(f_\alpha, A_E): (E,\alpha)\in [-4,4]\times [0,1]\}$,
the set of parameters corresponding to dynamics that are {\em not} uniformly hyperbolic.

Returning to spectral theory, 
we continue to explore the relationship between spectrum and dynamics.
If $\alpha$ is irrational, then $f_\alpha$ is ergodic, and Oseledets's theorem implies that the Lyapunov exponents for any cocycle over $f_\alpha$ take constant values over a full measure set.  Thus the Lyapunov exponents of $A_E$ over $f_\alpha$  take two essential values, $\chi^+_E\geq 0$, and $\chi^-_E$; the fact that $\det (A_E) = 1$ implies that $\chi^-_E = -\chi_E^+  \leq 0$.  Then either  $A_E$ is nonuniformly hyperbolic  (if $\chi^+_E>0$), or the exponents of $A_E$ vanish.

Thus for fixed $\alpha$ irrational, the spectrum  $\Sigma_{\alpha}$ splits, from a dynamical point of view, into two (measurable) sets: the set of $E$ for which $A_E$ is nonununiformly hyperbolic, and the set of $E$ for which the exponents of $A_E$ vanish.   On the other hand, 
spectral analysis gives us a different decomposition of the spectrum:
$$
\sigma(H^\alpha_x)=\sigma_{ac}(H^\alpha_x)\cup\sigma_{sc}(H^\alpha_x)\cup \sigma_{pp}(H^\alpha_x)
$$ 
where $\sigma_{ac}(H^\alpha_x)$ is the absolutely continuous spectrum, $\sigma_{pp}(H_x)$ is the pure point spectrum (i.e., the closure of the eigenvalues), and $\sigma_{sc}(H^\alpha_x)$ is the singular continuous spectrum.  All three types of spectra have meaningful physical interpretations.  While the spectrum $\sigma(H^\alpha_x)$ does not depend in $x$ (since $\alpha$ is irrational), the decomposition into subspectra can depend on $x$.\footnote{In fact, the decomposition  is independent of a.e. $x$, just not all $x$.}  It turns out that the absolutely continuous spectrum does not depend on $x$, so we can write $\Sigma_{ac,\alpha}$ for this common set.

The next deep relation between spectral theory and Lyapunov exponents is the following, which is due to Kotani:
\begin{theorem}[Kotani, \cite{kotani}]\label{t.kotani} Fix $\alpha$ irrational. Let $\cZ$ be the set of $E$ such that the Lyapunov exponents of $A_E$ over $f_\alpha$ vanish.
Let $\overline{\mathcal Z^{ess}}$ denote the essential closure of $\mathcal Z$, i.e. the closure of the Lebesgue density points of $\mathcal Z$. Then 
$$
\Sigma_{ac}=\overline{\mathcal Z^{ess}}.
$$
\end{theorem}

Thus Lyapunov exponents of the cocycle are  closely related to the spectral type of the operators $H_x$. For instance, Theorem~\ref{t.kotani} implies that if $A_E$ is nonuniformly hyperbolic over $f_\alpha$ for almost every $E\in \Sigma_\alpha$, then $\Sigma_{ac,\alpha}$ is empty:
$H^\alpha_x$ has no absolutely continuous spectrum.

We remark that Theorems~\ref{t.johnson} and \ref{t.kotani} hold  for much broader classes of 
Schr\"odinger operators over ergodic measure preserving systems.   For a short and self-contained proof of Theorem~\ref{t.johnson}, see \cite{zhang}.
The spectral theory of one-dimensional  Schr\"odinger operators is a rich subject, and we've only scratched the surface here; for further reading, see the recent surveys \cite{JM} and \cite{Dam}.

Avila's very recent work, some of it still unpublished, provides further
fascinating connections of this type, linking the spectral properties of  quasiperiodic operators with analytic
potentials to  properties of the
Lyapunov exponents of their associated cycles \cite{AGlobalActa, A0}.

\section{Spaces of dynamical systems and metadynamics}

Sections 2, 3 and 4  are all about {\em families} of dynamical systems.  In Section 2, the family is the space of all volume preserving diffeomorphisms of a compact manifold $M$.  This is an infinite dimensional, non-locally compact space, and we have thrown up our hands and depicted it in Figure 8 as a blob.  Theorem~\ref{t=acw} asserts that within a residual set of positive entropy systems (which turn out to be an  open subset of the  blob), measurable hyperbolicity (and ergodicity) is generic.

In Section 3, the moduli space $\cM$ of translation surfaces can also be viewed as a space of dynamical systems,  in particular the {\em billiard} flows on {\em rational} polygons,  i.e., polygons whose corner angles are multiples of $2\pi$.  In a billiard system, one shoots a ball in a fixed direction and records the location of the bounces on the walls.  By a process called unfolding, a billiard trajectory can be turned into a straight ray  in a translation surface.\footnote{Not every translation surface comes from a  billiard, since the billiards have extra symmetries. But the space of billiards embeds inside the space of translation surfaces, and the Teichm\"uller flow preserves the set of billiards.}   The process is illustrated in Figure 14 for the square torus billiard.
\begin{figure}[ht]
\begin{center}
\includegraphics[scale=0.4]{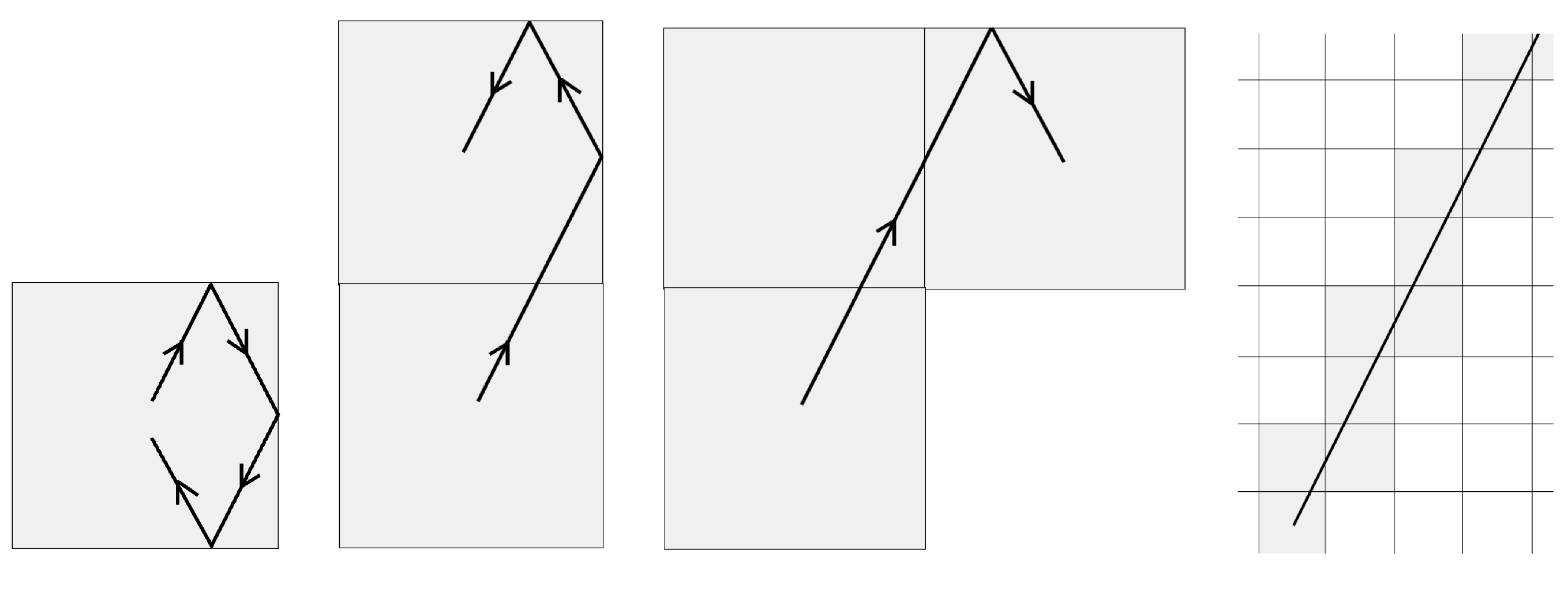}
\end{center}
\caption{Unfolding billiards in a square to get lines in a torus (courtesy Diana Davis).}
\end{figure}

The moduli space $\cM$ is not so easy to draw and not completely understood (except for $g=1$).  It is, however, a finite dimensional manifold and carries some nice structures, which makes it easier to
picture than $\Diff(M)$. Theorem~\ref{t=favc} illustrates how dynamical properties of a {\em meta dynamical system}, i.e. the Teichm\"uller flow $\cF_t\colon \cM\to \cM$, are tied to 
the dynamical properties of the elements of $\cM$.  For example, the Lyapunov exponents
of the KZ cocycle over $\cF_t$ for a given billiard table with a given direction describe how well
an infinite billiard ray can be approximated by closed, nearby billiard paths.

In Section 4, we saw how the spectral properties of a family of operators $\{H^\alpha_x: \alpha\in[0,1]\}$
 are reflected in the dynamical properties of families of
cocycles $\{(f_\alpha, A_E): (E,\alpha)\in [-4,4]\times [0,1]\}$.  Theorems about spectral properties thus have their dynamical counterparts.  For example,  Theorem~\ref{t.johnson}
tells us that the butterfly is the complement of those parameter values where the cocycle
$(f_\alpha, A_E)$ is uniformly hyperbolic. Since uniform hyperbolicity is an open property in both $\alpha$ and $E$, the complement of the butterfly is open.
  Corollary~\ref{c=meas0}  tells us that the butterfly has measure $0$.  Thus the set of parameter values in the square that are hyperbolic form an open and dense, full-measure subset.  
In fact, work of Bourgain-Jitomirskaya \cite{BJ} implies that  {\em the butterfly is precisely the set of parameter values $(E,\alpha)$ where the Lyapunov exponents of $(f_\alpha, A_E)$ vanish for some $x$.}\footnote{which
automatically means for all $x$ in case of irrational $\alpha$.}
These results in some ways  echo Theorem~\ref{t=acw}, within a very special family of dynamics.

The Hofstadter butterfly is just one instance of a low-dimensional family of dynamical systems containing interesting dynamics and rich structure.   
A similar picture appears in complex dynamics,\footnote{Another field in which Avila has made significant contributions, which we do not touch upon here.} in the $1$ (complex) parameter family of dynamical systems $\{p_c(z) = z^2 + c: c\in \CC\}$.  The Mandelbrot set consists of parameters $c$ for which the map $p_c$ has a connected Julia set $J_c$:
\begin{figure}[h]
\begin{center}
\includegraphics[scale=.42]{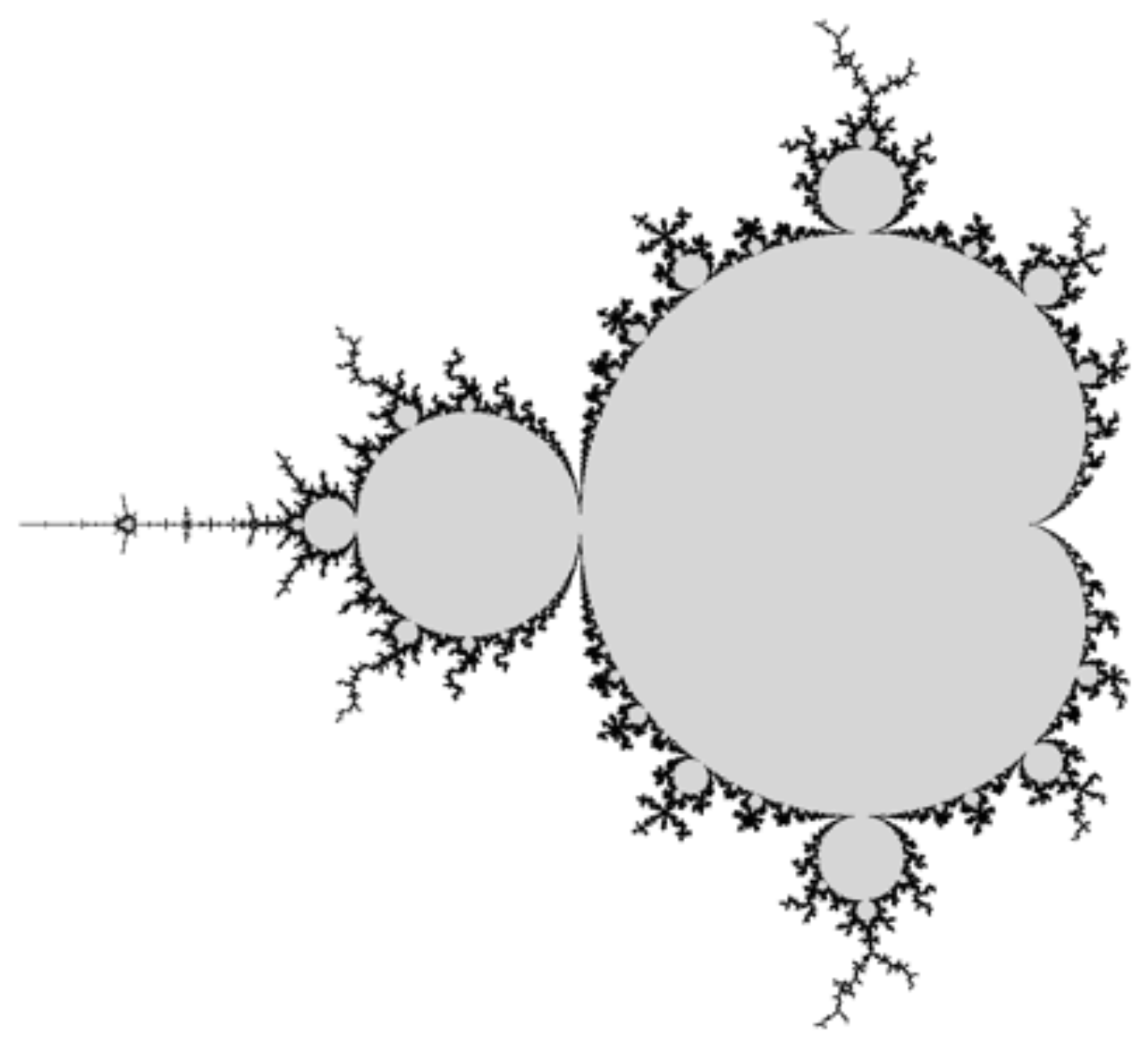}
\end{center}
\caption{The Mandelbrot Set}
\end{figure}

 Note that in this conformal context, uniform hyperbolicity of the derivative cocycle  of $p_c$ on $J_c$ just means that there exists an $n$ such that $|(p_c^n(z))'| \geq 2$, for all $z\in J_c$.
It is conjectured that the set of parameters $c$ such that $p_c$ is uniformly hyperbolic on $J_c$
is (open and) dense in the Mandelbrot set.

\section{Themes}

We conclude by summarizing a few themes, some of which  have come up in our discussion.

\bigskip
\noindent
{\bf Nonvanishing exponents sometimes produce chaotic behavior.}  The bedrock result in this regard is Anosov's proof \cite{Anosov}  that smooth Anosov flows and diffeomorphisms are mixing (and in particular ergodic).  
Another notable result is  Anatole Katok's proof \cite{katok} that measurable hyperbolicity of diffeomorphism  with respect to some measure $\mu$ produces many periodic orbits --- in particular, the number of orbits of period $n$ grows exponentially in $n$.  Pesin theory provides a sophisticated tool for exploiting measuable hyperbolicity to produce chaotic behavior such as mixing and even the Bernoulli property.

\bigskip

\noindent
{\bf  Exponents can carry geometric information.} We have not discussed it here, but there are delicate relationships between entropy, exponents and Hausdorff dimension of invariant sets and measures, established in full generality by Ledrappier-Young \cite{LY1, LY2}.
The  expository article \cite{lsy} contains a  clear discussion of these relationships as well as some of the other themes discussed in this paper.  The interplay between dimension, entropy and exponents has been fruitfully exploited in numerous contexts, notably in rigidity theory.  Some examples are:
Ratner's theorem for unipotent flows, Elon Lindenstrauss's proof of Quantum Unique Ergodicity for arithmetic surfaces, and the Einsiedler-Katok-Lindenstrauss proof that the set of exceptions to the Littlewood conjecture has Hausdorff dimension 0.  See \cite{Witte, ElonQUE, EKL}.

\bigskip

\noindent
{\bf Vanishing exponents sometimes present an exceptional situation that can be exploited.}
Both Furstenberg's theorem and Kotani theory illustrate this phenomenon.  Here's Furstenberg's criterion,
presented in a special case:
\begin{theorem}[Furstenberg, \cite{F}]\label{t=furstenberg} Let $(A_1,\ldots, A_k)\subset SL(2,\RR)$, and let $G$ be the smallest closed subgroup of $SL(2,\RR)$ containing $\{A_1,\ldots, A_k\}$.  Assume that:
\begin{enumerate}
\item $G$ is not compact.
\item  There is no finite collection of lines $\emptyset \neq  L \subset  \RR^2$
such that $M(L) = L$, for all $M \in G$.
\end{enumerate}
Then  for any probability vector $p=(p_1,\ldots p_k)$ on $\{1,\ldots, k\}$ with $p_i>0$, for all $i$, there exists $\chi^+(p)>0$, such that for  almost every $\omega\in \{1,\ldots, k\}^\NN$ (with respect to the Bernoulli measure
 $p^\NN$):
$$
\lim_{n\to\infty}\frac{1}{n}\log\|A^{(n)}(\omega)\| = \chi^+.
$$
\end{theorem}
One way to view this result:  if the exponent $\chi_+$ vanishes, then the matrices either leave invariant a common line or pair of lines, or they generate a precompact group.
Both possibilities are degenerate and are easily destroyed by perturbing the matrices.
One proof of a generalization of this result \cite{L} exploits the connections between entropy, dimension and exponents alluded to before.  This result was formulated more completely in a dynamical setting 
by \cite{BGV} as an ``Invariance Principle," which has been further refined and applied in various works of Avila and others.  See for example \cite{AV3, ASV, AVW}.  

For general $SL(d,\RR)$ cocycles, vanishing of exponents is still an exceptional situation, but even more generally, the condition $d_i \geq 2$ ---  that an exponent has multiplicity greater than 1 --- is also exceptional.   This statement was rigorously established for random matrix products by Guivarc'h-Raugi \cite{GR} and undleries 
the Avila-Viana proof of simplicity of spectrum for the KZ cocycle.  See the discussion at the end of Section 3.
The same ideas  play a role in Margulis's proof of superrigidity for higher rank lattices in semisimple Lie groups.
See \cite{Msuper}.

\bigskip

\noindent
{\bf Continuity and regularity of exponents is a tricky business.} 
In general, Lyapunov exponents do not depend smoothly, or even continuously, on the cocycle.  Understanding the 
exact relationship between exponents, cocycles, measures and dynamics is an area still under exploration, 
and a few of Avila's deepest results, some of them still in preparation with Eskin and Viana,  lie in this area. 
The book \cite{Vianabook} is an excellent introduction to the subject.
\bigskip

\noindent{\bf Acknowledgments.} Effusive thanks to  Artur Avila, 
 Svetlana Jitomirskaya, Curtis McMullen,  Zhenghe Zhang,  and Anton Zorich for patiently explaining a lot of math to me, to  Clark Butler, Kathryn Lindsay, Kiho Park, Jinxin Xue and Yun Yang for catching many errors, and to Diana Davis, Carlos Matheus, Curtis McMullen and Marcelo Viana for generously sharing their images.
I am also indebted to Bryna Kra for carefully reading a draft of this paper and suggesting numerous improvements.

\end{document}